\input{epsf}

\newcounter{letter}     
\newenvironment{alphalist}{\begin{list}
{{\normalshape(\alph{letter})}}{\usecounter{letter}}}{\end{list}}

        \newcommand{\be}{\begin{equation}}
        \newcommand{\ee}{\end{equation}}
        \newcommand{\ba}{\begin{eqnarray}}
        \newcommand{\ea}{\end{eqnarray}}
        \newcommand{\ban}{\begin{eqnarray*}}
        \newcommand{\ean}{\end{eqnarray*}}
        \newcommand{\barr}{\begin{array}}
        \newcommand{\earr}{\end{array}}

\def\=>{\Rightarrow}
\newcommand{\To}{\Rightarrow}
\def\mapright#1{\smash{\mathop{\longrightarrow}\limits^{#1}}}

\newcommand{\et}{\hspace{-0.08in}{\bf .}\hspace{0.1in}}
\newcommand{\BOX}{\hbox {$\sqcap$ \kern -1em $\sqcup$}}

\renewcommand{\hom}{{\rm hom}}

\newcommand{\FinSet}{{\rm FinSet}}
\newcommand{\Gpd}{{\rm Gpd}}
\newcommand{\Top}{{\rm Top}}
\newcommand{\Type}{{\rm Type}}
\newcommand{\Cat}{{\rm Cat}}

\newcommand{\Cob}{{\rm Cob}}
\newcommand{\Tang}{{\rm Tang}}
\newcommand{\Braid}{{\rm Braid}}
\renewcommand{\to}{\rightarrow}
\newcommand{\tensor}{\otimes}
\newcommand{\maps}{\colon}

\newcommand{\iso}{\cong}
\newcommand{\elt}{{\rm elt}}

\newcommand{\R}{{\Bbb R}}
\newcommand{\N}{{\Bbb N}}
\newcommand{\Z}{{\Bbb Z}}

\newcommand{\Aut}{{\rm aut}}
\newcommand{\End}{{\rm end}}
\newcommand{\Decat}{{\rm Decat}}
\newcommand{\Disc}{{\rm Disc}}

\newtheorem{thm}{Theorem}    

\newtheorem{defn}[thm]{Definition}

	\documentstyle[12pt,diagram,amsfonts]{article}
	\begin{document}

      \begin{center}
      {\bf Categorification \\}
      \vspace{0.5cm}
      {\em John C.\ Baez \\}
      \vspace{0.3cm}
      {\small Department of Mathematics, University of California\\
      Riverside, California 92521 \\
      USA\\ }
      \vspace{0.3cm}
      {\em James Dolan \\}
      \vspace{0.3cm}
      {\small Sanga Research Corporation \\
      2015 Rue Peel, Suite 1000 \\
      Montreal, Quebec H3A 1T8 \\
      Canada \\ }
      \vspace{0.3cm}
      {\small email: baez@math.ucr.edu, jdolan@sangacorp.com \\}
      \vspace{0.3cm}
      {\small February 2, 1998 \\ }
      \end{center}

\begin{abstract}

Categorification is the process of finding category-theoretic analogs of
set-theoretic concepts by replacing sets with categories, functions with
functors, and equations between functions by natural isomorphisms
between functors, which in turn should satisfy certain equations of
their own, called `coherence laws'.   Iterating this process requires a
theory of `$n$-categories', algebraic structures having objects,
morphisms between objects, 2-morphisms between morphisms and so on up to
$n$-morphisms.  After a brief introduction to $n$-categories and their
relation to homotopy theory, we discuss algebraic structures that can be
seen as iterated categorifications of the natural numbers and integers.
These include tangle $n$-categories, cobordism $n$-categories, and
the homotopy $n$-types of the loop spaces $\Omega^k S^k$.  We conclude
by describing a definition of weak $n$-categories based on the theory
of operads.  
\end{abstract}

\section{Introduction}

The tongue-twisting term `categorification' was invented by Crane
\cite{Crane,CF}.  It refers to the process of finding category-theoretic
analogs of ideas phrased in the language of set theory, using the
following analogy between set theory and category theory:

\vskip 2em
\begin{center}
{\small
\begin{tabular}{|c|c|}                    \hline
elements    &  objects                       \\     \hline
equations   &  isomorphisms                  \\     
between elements   &  between objects        \\     \hline
sets        &  categories                    \\     \hline
functions   &  functors                      \\     \hline
equations   &  natural isomorphisms    \\     
between functions   & between functors       \\     \hline
\end{tabular}} \vskip 1em
1.  Analogy between set theory and category theory
\end{center}
\vskip 0.5em

\noindent Just as sets have elements, categories have objects.  Just as
there are functions between sets, there are functors between categories.
Interestingly, the proper analog of an equation between elements is not
an equation between objects, but an isomorphism.   More generally, the
analog of an equation between functions is a natural isomorphism between
functors.  

For example, the category $\FinSet$, whose objects are finite sets and
whose morphisms are functions, is a categorification of the set $\N$ of
natural numbers.   The disjoint union and Cartesian product of finite
sets correspond to the sum and product in $\N$, respectively.  Note that
while addition and multiplication in $\N$ satisfy various equational
laws such as commutativity, associativity and distributivity, disjoint
union and Cartesian product satisfy such laws {\it only up to natural
isomorphism}.  

If one studies categorification one soon discovers an amazing fact: many
deep-sounding results in mathematics are just categorifications of facts
we learned in high school!  There is a good reason for this.  All along,
we have been unwittingly `decategorifying' mathematics by pretending
that categories are just sets.  We `decategorify' a category by
forgetting about the morphisms and pretending that isomorphic objects
are equal.  We are left with a mere set: the set of isomorphism classes
of objects. 

To understand this, the following parable may be useful.  Long ago, when
shepherds wanted to see if two herds of sheep were isomorphic, they
would look for an explicit isomorphism.  In other words, they would line
up both herds and try to match each sheep in one herd with a sheep in
the other.  But one day, along came a shepherd who invented
decategorification.  She realized one could take each herd and `count'
it, setting up an isomorphism between it and some set of `numbers',
which were nonsense words like `one, two, three, \dots' specially
designed for this purpose.  By comparing the resulting numbers, she
could show that two herds were isomorphic without explicitly
establishing an isomorphism!  In short, by decategorifying the category
of finite sets, the set of natural numbers was invented.   

According to this parable, decategorification started out as a stroke of
mathematical genius.  Only later did it become a matter of dumb habit,
which we are now struggling to overcome by means of categorification.
While the historical reality is far more complicated, categorification
really has led to tremendous progress in mathematics during the 20th
century.  For example, Noether revolutionized algebraic topology by
emphasizing the importance of homology groups.  Previous work had
focused on Betti numbers, which are just the dimensions of the rational
homology groups.  As with taking the cardinality of a set, taking the
dimension of a vector space is a process of decategorification, since
two vector spaces are isomorphic if and only if they have the same
dimension.  Noether noted that if we work with homology groups rather
than Betti numbers, we can solve more problems, because we obtain
invariants not only of spaces, but also of maps.  In modern parlance,
the $n$th rational homology is a {\it functor} defined on the {\it
category} of topological spaces, while the $n$th Betti number is a mere
{\it function} defined on the {\it set} of isomorphism classes of
topological spaces.  Of course, this way of stating Noether's insight is
anachronistic, since it came before category theory.  Indeed, it was in
Eilenberg and Mac Lane's subsequent work on homology that category
theory was born!

Decategorification is a straightforward process which typically 
destroys information about the situation at hand.  Categorification,
being an attempt to recover this lost information, is inevitably fraught
with difficulties.  One reason is that when categorifying, one does not
merely replace equations by isomorphisms.  One also demands that these
isomorphisms satisfy some new equations of their own, called `coherence
laws'.  Finding the right coherence laws for a given situation is
perhaps the trickiest aspect of categorification.

For example, a monoid is a set with a product satisfying the associative
law and a unit element satisfying the left and right unit laws.   The
categorified version of a monoid is a `monoidal category'.   This is a
category $C$ with a product $\tensor \maps C \times C \to C$ and and
unit object $1 \in C$.   If we naively impose associativity and the
left and right unit laws as equational laws, we obtain the definition
of a `strict' monoidal category.   However, the philosophy of
categorification  suggests instead that we impose them only up to
natural isomorphism.   Thus, as part of the structure of a `weak' monoidal
category, we specify a natural isomorphism
\[   a_{x,y,z}\maps (x \tensor y) \tensor z \to x \tensor (y \tensor z)  \]
called the `associator', together with natural isomorphisms
\[                l_x \maps 1 \tensor x \to x, \]
\[                 r_x \maps x \tensor 1 \to x .\]
Using the associator one can construct isomorphisms between any two
parenthesized versions of the tensor product of several objects.  However,
we really want a {\it unique} isomorphism.  For example, there are 5
ways to parenthesize the tensor product of 4 objects, which 
are related by the associator as follows:
\[    \begin{diagram}[((x \tensor y)\tensor z)\tensor w]
\node{((x \tensor y) \tensor z)\tensor w}
\arrow{e,t}{a_{x\tensor y,z,w}}
\arrow{s,l}{a_{x,y,z}\tensor w}
\node{(x\tensor y)\tensor(z\tensor w)}
\arrow{e,t}{a_{x,y,z\tensor w}}
\node{x \tensor (y\tensor(z \tensor w))}  \\
\node{(x \tensor (y \tensor z)) \tensor w}
\arrow[2]{e,t}{a_{x,y\tensor z,w}}
\node[2]{x \tensor ((y \tensor z)\tensor w)}
\arrow{n,r}{x \tensor a_{y,z,w}}
\end{diagram} \]
In the definition of a weak monoidal category we impose a coherence
law, called the `pentagon identity', saying that this diagram commutes.  
Similarly, we impose a coherence law saying that the following diagram
built using $a,l$ and $r$ commutes:
\[      \begin{diagram}[(1 \tensor x) \tensor 1]
\node{(1 \tensor x) \tensor 1} \arrow{s,l}{l_x \tensor 1} 
\arrow[2]{e,t}{a_{1,x,1}} \node[2]{1 \tensor (x\tensor 1)} \arrow{s,r}{1
\tensor r_x} \\
\node{x \tensor 1} \arrow{e,t}{r_x} \node{x} \node{1 \tensor x}
\arrow{w,t}{l_x} 
\end{diagram} \]

This definition raises an obvious question: how do we know we have
found all the right coherence laws?  Indeed, what does `right'
even {\it mean} in this context?  Mac Lane's coherence theorem
\cite{MacLane} gives one answer to this question: the above coherence
laws imply that any two isomorphisms built using $a,$ $l$ and $r$
and having the same source and target must be equal. 

Further work along these lines allow us to make more precise the sense
in which $\N$ is a decategorification of $\FinSet$.  For example, just
as $\N$ forms a monoid under either addition or multiplication,
$\FinSet$ becomes a monoidal category under either disjoint union or
Cartesian product if we choose the isomorphisms $a,l$, and $r$
sensibly.  In fact, just as $\N$ is a `rig', satisfying all the ring
axioms except those involving additive inverses, $\FinSet$ is what one
might call a `rig category'.  In other words, it satisfies the rig
axioms up to natural isomorphisms satisfying the coherence laws
discovered by Kelly \cite{Kelly} and Laplaza \cite{Laplaza}, who proved
a coherence theorem in this context.  Just as the decategorification of
a monoidal category is a monoid, the decategorification of any rig
category is a rig.  In particular, decategorifying the rig category
$\FinSet$ gives the rig $\N$.  This idea is especially important in
combinatorics, where the best proof of an identity involving natural numbers
is often a `bijective proof': one that actually establishes an isomorphism
between finite sets \cite{Joyal,Rota}.

While coherence laws can sometimes be justified retrospectively by
coherence theorems, certain puzzles point to the need for a deeper
understanding of the {\it origin} of coherence laws.   For example,
suppose we want to categorify the notion of `commutative monoid'.   The
strictest possible approach, where we take a strict monoidal category
and impose an equational law of the form $x \tensor y = y \tensor x$, is
almost completely uninteresting.  It is much better to start with a weak
monoidal category equipped with a natural isomorphism $B_{x,y} \maps x
\tensor y \to y \tensor x$ called the `braiding', and then impose 
coherence laws called `hexagon identities' saying that the following 
two diagrams commute:
\[ \begin{diagram}[(x \tensor y) \tensor z)]
\node{x \tensor (y \tensor z)} \arrow{s,l}{B_{x,y \tensor z}}
\arrow{e,t}{a^{-1}_{x,y,z}}
\node{(x \tensor y) \tensor z} \arrow{e,t}{B_{x,y} \tensor z}
\node{(y \tensor x) \tensor z} \arrow{s,r}{a_{y,x,z}} \\
\node{(y \tensor z) \tensor x} \arrow{e,b}{a_{y,z,x}}
\node{y \tensor (z \tensor x)} \arrow{e,t}{y \tensor B_{x,z}} 
\node{y \tensor (x \tensor z)} 
\end{diagram} \]
\[ \begin{diagram}[(x \tensor y) \tensor z)]
\node{(x \tensor y) \tensor z} \arrow{s,l}{B_{x \tensor y, z}}
\arrow{r,t}{a_{x,y,z}} 
\node{x \tensor (y \tensor z)} \arrow{e,t}{x \tensor B_{y,z}}
\node{x \tensor (z \tensor y)} \arrow{s,r}{a_{x,z,y}^{-1}} \\
\node{z \tensor (x \tensor y)} \arrow{e,b}{a^{-1}_{z,x,y}}
\node{(z \tensor x) \tensor y}  \arrow{e,t}{B_{x,z} \tensor y} 
\node{(x \tensor z) \tensor y}
\end{diagram} \]
This gives the definition of a weak `braided monoidal category'.  If we
impose an additional coherence law of the form $B_{x,y}^{-1} = B_{y,x}$,
we obtain the definition of a `symmetric monoidal category'.  
Both of these concepts are very important; which one is `right' depends
on the context.   However, neither implies that every pair of parallel
morphisms built using the braiding are equal.   A good theory of
coherence laws must naturally account for these facts.

The deepest insights into such puzzles have traditionally come from
topology.  In homotopy theory it causes problems to work with spaces
equipped with algebraic structures satisfying equational laws, because
one cannot transport such structures along homotopy equivalences.   It
is better to impose laws {\it only up to homotopy}, with these
homotopies satisfying certain coherence laws, but again only up to
homotopy, with these higher homotopies satisfying their own higher
coherence laws, and so on.  Coherence laws thus arise naturally in
infinite sequences.  For example, Stasheff \cite{Stasheff} discovered
the pentagon identity and a sequence of higher coherence laws for
associativity when studying the algebraic structure possessed by a space
that is homotopy equivalent to a loop space.  Similarly, the hexagon
identities arise as part of a sequence of coherence laws for spaces
homotopy equivalent to double loop spaces, while the extra coherence law
for symmetric monoidal categories arises as part of a sequence for
spaces homotopy equivalent to triple loop spaces.   The higher coherence
laws in these sequences turn out to be crucial when we try to {\it
iterate} the process of categorification.  

Starting in the late 1960's, Boardman, Vogt \cite{BV,BV2} and others
developed the study of these higher coherence laws into a full-fledged
theory of `homotopy-invariant algebraic structures'.   However, we have
yet to attain a general and systematic  theory of categorification,
particularly when it comes to iterated categorification, which requires
a good theory of `$n$-categories'.  The main goal of this paper is to
outline our current understanding of categorification and to point out
some directions for further study.  

The plan of the paper is as follows.  In Section \ref{ncats} we give a
quick tour of $n$-category theory.  This theory is just beginning to be
developed, and there are various alternative approaches which have not
yet been reconciled, but here we leave out most of the technical details
and sketch what we expect from any reasonable approach.  In Section
\ref{lessons} we discuss in more detail the lessons homotopy theory has
for $n$-category theory.  We hope this section can be followed even 
by those who are not already experts on homotopy theory.   
In Section \ref{examples} we describe some algebraic structures that 
amount to iterated categorifications of the natural numbers and the integers.  
A large amount of interesting mathematics emerges from the study of these
structures.  In Section \ref{opetopes} we summarize our own approach to
$n$-categories.  

In a previous paper \cite{BD} we sketched a program of using
$n$-categories to clarify the relationships between topological quantum
field theory and more traditional approaches to algebraic topology.  The
present paper covers some aspects of this program in more detail, taking
advantage of work that has been done in the meantime.  Various other
aspects are treated in a series of papers entitled `Higher-Dimensional
Algebra' \cite{B,BD2,BL,BN}.

\section{$n$-Categories} \label{ncats}

One philosophical reason for categorification is that it refines our
concept of `sameness' by allowing us to distinguish between isomorphism
and equality.  In a set, two elements are either the same or different. 
In a category, two objects can be `the same in a way' while still being
different.  In other words, they can be isomorphic but not equal.  Even
more importantly, two objects can be the same in more than one way,
since there can be different isomorphisms between them.   This gives
rise to the notion of the `symmetry group' of an object: its group of
automorphisms. 

Consider, for example, the fundamental groupoid $\Pi_1(X)$ of a 
topological space $X$: the category with points of $X$ as objects and
homotopy classes of paths with fixed endpoints as morphisms.  This
category captures all the homotopy-theoretic information about $X$ in
dimensions $\le 1$ --- or more precisely, its homotopy 1-type
\cite{Baues}.   The group of automorphisms of an object $x$ in this
category is just the fundamental group $\pi_1(X,x)$.  If we decategorify
the fundamental groupoid of $X$, we forget {\it how} points in $X$ are
connected by paths, remembering only {\it whether} they are, and we
obtain the set of components of $X$.   This captures only the homotopy
0-type of $X$.  

This example shows how decategorification eliminates `higher-dimensional
information' about a situation.   Categorification is an attempt to
recover this information.  This example also suggests that we can keep
track of the homotopy 2-type of $X$ if we categorify further and 
distinguish between paths that are equal and paths that are merely
isomorphic (i.e., homotopic).   For this we should work with a
`2-category' having points of $X$ as objects, paths as morphisms, and
certain equivalence classes of homotopies between paths as 2-morphisms. 

In a marvelous self-referential twist, the definition of `2-category' is
simply the categorification of the definition of `category'!   Like a
category, a 2-category has a class of objects, but now for any pair
$x,y$ of objects there is no longer a set $\hom(x,y)$; instead, there is
a category $\hom(x,y)$.  Objects of $\hom(x,y)$ are called morphisms of
$C$, and morphisms between them are called 2-morphisms of $C$. 
Composition is no longer a function, but rather a functor:
\[             \circ \maps \hom(x,y) \times \hom(y,z) \to \hom(x,z) .\]
For any object $x$ there is an identity $1_x \in \hom(x,x)$.  And now we
have a choice.  On the one hand, we can impose associativity and the
left and right unit laws strictly, as equational laws.  If we do this,
we obtain the definition of `strict 2-category' \cite{KS}.  On the other
hand, we can impose them only up to natural isomorphism, with these
natural isomorphisms satisfying the coherence laws discussed in the
previous section.  This is clearly more compatible with the spirit of
categorification.  If we do this, we obtain the definition of `weak
2-category' \cite{Benabou}.  (We warn the reader that strict
2-categories are traditionally known as `2-categories', while weak
2-categories are known as `bicategories'.  The present style of
terminology, introduced by Kapranov and Voevodsky \cite{KV}, has the
advantage of generalizing easily to $n$-categories for arbitrary $n$.)

The classic example of a 2-category is $\Cat$, which has categories as
objects, functors as morphisms, and natural transformations as
2-morphisms.   The presence of 2-morphisms gives $\Cat$ much of its
distinctive flavor, which we would miss if we treated it as a mere
category.  Indeed, Mac Lane has said that categories were originally
invented, not to study functors, but to study natural transformations! 
A good example of two functors that are not equal, but only naturally
isomorphic, are the identity functor and the `double dual' functor on
the category of finite-dimensional vector spaces.  

Given a topological space $X$, we can form a 2-category $\Pi_2(X)$
called the `fundamental 2-groupoid' of $X$.  The objects of this 
2-category are the points of $X$.   Given $x,y \in X$, the morphisms
from $x$ to $y$ are the paths $f \maps [0,1] \to X$ starting at $x$ and
ending at $y$.  Finally, given $f,g \in \hom(x,y)$, the 2-morphisms from
$f$ to $g$ are the homotopy classes of paths in $\hom(x,y)$ starting at
$f$ and ending at $g$.    Since the associative law for composition of
paths holds only up to homotopy, this 2-category is a weak 2-category. 
One can prove that it captures the homotopy 2-type of $X$ \cite{BS,MW}.
If we decategorify the fundamental 2-groupoid of $X$, we obtain its
fundamental groupoid.

From 2-categories it is a short step to dreaming of $n$-categories and
even $\omega$-categories --- but it is not so easy to make these dreams
into smoothly functioning mathematical tools.   Roughly speaking, an
$n$-category should be some sort of algebraic structure having objects,
1-morphisms between objects, 2-morphisms between 1-morphisms, and so on
up to $n$-morphisms.  There should be various ways of composing
$j$-morphisms for $1 \le j \le n$, and these should satisfy various
laws.  As with 2-categories, we can try to impose these laws either
strictly or weakly. Strict $n$-categories have been understood for quite
some time now \cite{Crans,EK}, but more interesting for us are the weak
ones.   Various definitions of weak $n$-category are currently under active 
study \cite{BD2,Batanin,HMP,Simpson,Simpson2,Street,Tamsamani,Tamsamani2}, 
and we discuss our own in Section \ref{opetopes}.   Here, however, we
wish to sketch the main challenges any theory of weak $n$-categories
must face, and some of the richness inherent in the notion of weak
$n$-category.    

{\it Nota bene: Throughout the rest of this paper, `$n$-category' will
mean `weak $n$-category' unless otherwise specified, and similarly for
`$n$-groupoid', etc..}   For the sake of definiteness, we shall
temporarily speak in terms of the `globular' approach to $n$-categories
developed by Batanin \cite{Batanin}.  In this approach, for $j \ge 1$
any $j$-morphism $a \maps f \to g$ has a source $f$ and a target $g$,
which are $(j-1)$-morphisms.  When $j \ge 2$ we require that $f$ and $g$
are parallel, by which we mean that $f,g \maps x \to y$ for some
$(j-2)$-morphisms $x$ and $y$.   In this approach, we visualize a
$j$-morphism as a $j$-dimensional disc.   For example, for $j = 2$:

\medskip
\centerline{\epsfysize=1.0in\epsfbox{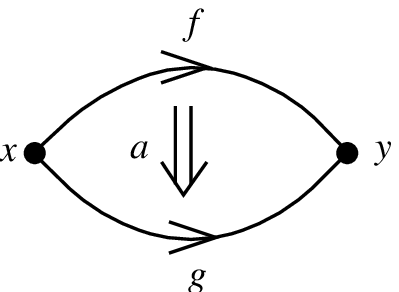}}
\medskip

\noindent Other approaches to $n$-categories use $j$-morphisms with
other shapes, such as simplices, discussed in Section \ref{lessons}, or
opetopes, discussed in Section \ref{opetopes}.  We believe that there is
basically a single notion of weak $n$-category lurking behind these
different approaches.  If this is true, they will eventually be shown to
be equivalent, and choosing among them will be merely a matter of
convenience.  However, the precise meaning of `equivalence' here is
itself rather subtle and $n$-categorical in flavor \cite{B2,BD2}.

The first challenge to any theory of $n$-categories is to give an
adequate treatment of coherence laws.  Composition in an $n$-category
should satisfy equational laws only at the top level, between
$n$-morphisms.   Any law concerning $j$-morphisms for $j < n$ should
hold only `up to equivalence'.   Here a $n$-morphism is defined to be an
`equivalence' if it is invertible, while for $j < n$ a $j$-morphism is
recursively defined to be an equivalence if it is invertible {\it up to
equivalence}.  Equivalence is generally the correct substitute for the
notion of equality in $n$-categorical mathematics.  When laws are
formulated as equivalences, these equivalences should in turn satisfy
coherence laws of their own, but again only up to equivalence, and so
on.  This becomes ever more complicated and unmanageable with increasing
$n$ unless one takes a systematic approach to coherence laws.  

The second challenge to any theory of $n$-categories is to handle
certain key examples.  First, for any $n$, there should be an
$(n+1)$-category $n\Cat$, whose objects are (small) $n$-categories, whose
morphisms are suitably weakened functors between these, whose
2-morphisms are suitably weakened natural transformations, and so on. 
Here by `suitably weakened' we refer to the fact that all laws should
hold only up to equivalence.  Second, for any topological space $X$,
there should be an $n$-category $\Pi_n(X)$ whose objects are points of
$X$, whose morphisms are paths, whose 2-morphisms are paths of paths,
and so on, where we take homotopy classes only at the top level. 
$\Pi_n(X)$ should be an `$n$-groupoid', meaning that all its
$j$-morphisms are equivalences for $0 \le j \le n$.  We call $\Pi_n(X)$
the `fundamental $n$-groupoid of $X$'.   Conversely, any $n$-groupoid
should determine a topological space, its `geometric realization'.  In
fact, these constructions should render the study of $n$-groupoids
equivalent to that of homotopy $n$-types (in a certain sense to be made
precise in the next section).  

A bit of the richness inherent in the concept of $n$-category becomes
apparent when we make the following observation: an $(n+1)$-category
with only one object can be regarded as special sort of $n$-category. 
Suppose that $C$ is an $(n+1)$-category with one object $x$.  Then we
can form the $n$-category $\tilde C$ by re-indexing: the objects of
$\tilde C$ are the morphisms of $C$, the morphisms of $\tilde C$ are the
2-morphisms of $C$, and so on.  The $n$-categories we obtain this way
have extra structure.  In particular, since the objects of $\tilde C$ are
really morphisms in $C$ from $x$ to itself, we can `multiply' (that is,
compose) them.

The simplest example is this: if $C$ is a category with a single object
$x$, $\tilde C$ is the set of endomorphisms of $x$.  This set is
actually a monoid.  Conversely, any monoid can be regarded as the monoid
of endomorphisms of $x$ for some category with one object $x$.  We
summarize this situation by saying that `a one-object category is a
monoid'.   Similarly, a one-object 2-category is a monoidal category. 
It is natural to expect this pattern to continue in all higher
dimensions; in fact, it is probably easiest to cheat and {\it define} a
monoidal $n$-category to be an $(n+1)$-category with one object.

Things get even more interesting when we iterate this process.  Given an
 $(n+k)$-category $C$ with only one object, one morphism, and so on up
to one $(k-1)$-morphism, we can form an $n$-category whose $j$-morphisms
are the $(j+k)$-morphisms of $C$.  In doing so we obtain a particular
sort of $n$-category with extra structure and properties, which we call
a `$k$-tuply monoidal' $n$-category.  Table 2 shows what we expect these
to be like for low values of $n$ and $k$.   For example, the
Eckmann-Hilton argument \cite{BD,BN,EH} shows that a 2-category with one
object and one morphism is a commutative monoid.  Categorifying this
argument, one can show that a 3-category with one object and one
morphism is a braided monoidal category.  Similarly, we expect that a
4-category with one object, one morphism and one 2-morphism is a
symmetric monoidal category, though this has not been worked out in full
detail, because of our poor understanding of 4-categories.   The fact
that both braided and symmetric monoidal categories appear in this table
seems to explain why both are natural concepts.  

\vskip 0.5em
\begin{center}
{\small
\begin{tabular}{|c|c|c|c|}  \hline
         & $n = 0$   & $n = 1$    & $n = 2$          \\     \hline
$k = 0$  & sets      & categories & 2-categories     \\     \hline
$k = 1$  & monoids   & monoidal   & monoidal         \\
         &           & categories & 2-categories     \\     \hline
$k = 2$  &commutative& braided    & braided          \\
         & monoids   & monoidal   & monoidal         \\
         &           & categories & 2-categories     \\     \hline
$k = 3$  &`'         & symmetric  & weakly involutory \\
         &           & monoidal   & monoidal         \\
         &           & categories & 2-categories     \\     \hline
$k = 4$  &`'         & `'         & strongly involutory\\
         &           &            & monoidal         \\
         &           &            & 2-categories     \\     \hline
$k = 5$  &`'         &`'          & `'               \\
         &           &            &                  \\
         &           &            &                  \\     \hline
\end{tabular}} \vskip 1em
2.  $k$-tuply monoidal $n$-categories
\end{center}
\vskip 0.5em

In any reasonable approach to $n$-categories there should be an
$n$-category $n\Cat_k$ whose objects are $k$-tuply monoidal weak
$n$-categories.  One should also be able to treat $n\Cat_k$ as a full
sub-$(n+k)$-category of $(n+k)\Cat$, though even for low $n,k$ this is
perhaps not as well known as it should be.  Consider for example $n = 0,
k = 1$.  The objects of $0\Cat_1$ are one-object categories, or monoids.
The morphisms of $0\Cat_1$ are functors between one-object categories,
or monoid homomorphisms.  But $0\Cat_1$ also has 2-morphisms 
corresponding to natural transformations.  We leave it as an exercise to
work out what these are in concrete terms.   More recently, Kapranov and
Voevodsky \cite{KV} have considered the case $n = k = 1$.  

We gave a detailed discussion of Table 2 in an earlier paper \cite{BD},
and subsequent work by various authors has improved our understanding of
some of the higher entries \cite{BN,Crans2,DS}.   There are many
interesting processes going from each entry in this table to its
neighbors.  We list some of the main ones below.  Most of these have
only been thoroughly studied for low values of $n$ and $k$, often in the
framework of `semistrict' $n$-categories, which are a kind of halfway
house between strict and weak ones.  We expect that they all generalize
to weak $k$-tuply monoidal $n$-categories for arbitrary $n$ and $k$, but
in many cases this has not yet been proved.   

\vskip 1em $\bullet$ \noindent
{\it Decategorification: $(n,k) \to (n-1,k)$}.  
Let $C$ be a $k$-tuply monoidal $n$-category $C$.  Then there
should be a $k$-tuply monoidal $(n-1)$-category $\Decat C$ whose
$j$-morphisms are the same as those of $C$ for $j < n - 1$, but whose
$(n-1)$-morphisms are isomorphism classes of $(n-1)$-morphisms of $C$.  

\vskip 1em $\bullet$ \noindent
{\it Discrete categorification: $(n,k) \to (n+1,k)$}.  
There should be a `discrete' $k$-tuply monoidal $(n+1)$-category
$\Disc C$  having the $j$-morphisms of $C$ as its $j$-morphisms for $j
\le n$, and only identity $(n+1)$-morphisms.  The decategorification of
$\Disc C$ should be $C$.

\vskip 1em $\bullet$ \noindent
{\it Delooping: $(n,k) \to (n+1,k-1)$}.   
There should be a $(k-1)$-tuply monoidal $(n+1)$-category $BC$ with
one object obtained by reindexing, the $j$-morphisms of $BC$ being the
$(j+1)$-morphisms of $C$.   We use the notation `$B$' and call $BC$ the
`delooping' of $C$ because of its relation to the classifying space 
construction in topology.   

\vskip 1em $\bullet$ \noindent
{\it Looping: $(n,k) \to (n-1,k+1)$}.  
Given objects $x,y$ in an $n$-category, there should be an
$(n-1)$-category $\hom(x,y)$.  If $x = y$ this should be a monoidal
$(n-1)$-category, and we denote it as $\End(x)$.  For $k > 0$, if $1$
denotes the unit object of the $k$-tuply monoidal $n$-category $C$,
$\End(1)$ should be a $(k+1)$-tuply monoidal $(n-1)$-category.  We call
this process `looping', and denote the result as $\Omega C$, because of
its relation to loop space construction in topology.  For $k > 0$,
looping should extend to an $(n+k)$-functor $\Omega\maps n\Cat_k \to
(n-1)\Cat_{k+1}$.   The case $k = 0$ is a bit different: we should be
able to loop a `pointed' $n$-category, one having a distinguished object
$x$, by letting $\Omega C = \End(x)$.   In either case, the
$j$-morphisms of $\Omega C$ correspond to certain $(j-1)$-morphisms of
$C$.

\vskip 1em $\bullet$ \noindent
{\it Forgetting monoidal structure: $(n,k) \to (n,k-1)$}.  
By forgetting the $k$th level of monoidal structure, we should be able
to think of $C$ as a $(k-1)$-tuply monoidal $n$-category $FC$.  This
should extend to an $n$-functor $F \maps n\Cat_k \to n\Cat_{k-1}$.

\vskip 1em $\bullet$ \noindent
{\it Stabilization: $(n,k) \to (n,k+1)$}. 
Though adjoint $n$-functors are still poorly understood, there should
be a left adjoint to forgetting monoidal structure, which we call
`stabilization' and denote by $S \maps n\Cat_k \to n\Cat_{k+1}$.  
(In our previous work we called it `suspension', but this is probably a
bit misleading.)   The Stabilization Hypothesis \cite{BD} states that
for $k \ge n + 2$, stabilization is an equivalence from $n\Cat_k$ to
$n\Cat_{k+1}$.   This is why the $n$th column of Table 2 has only $n+2$
distinct entries, and then settles down.  While not yet proven or even
formulated as a precise conjecture except in low dimensions, there is a
lot of good evidence for this hypothesis, some of which we mention in 
the next section.  In what follows, we assume this hypothesis and call a
$k$-tuply monoidal $n$-category with $k = n + 2$ a `stable $n$-category'.

\vskip 1em $\bullet$ \noindent
{\it Forming the generalized center: $(n,k) \to (n,k+1)$}.  
Thinking of $C$ as an object of the $(n+k)$-category $n\Cat_k$, there
should be a $(k+1)$-tuply monoidal $n$-category $ZC$, the `generalized
center' of $C$, given by $\Omega^k(\End(C))$.  In other words, $ZC$ is
the largest sub-$(n+k+1)$-category of $(n+k)\Cat$ having $C$ as its only
object, $1_C$ as its only morphism, $1_{1_C}$ as its only 2-morphism,
and so on up to dimension $k$.  This construction gets its name from the
case $n = 0$, $k = 1$, where $ZC$ is the usual center of the monoid
$C$.  Categorifying leads to the case $n = 1$, $k = 1$, which gives a very
important construction of braided monoidal categories from monoidal
categories \cite{JS,KV,Majid}.  In particular, when $C$ is the monoidal
category of representations of a Hopf algebra $H$, $ZC$ is the braided
monoidal category of representations of the quantum double $D(H)$. 
Categorifying still further, Baez and Neuchl \cite{BN} treated the case
$n = 2$, $k = 1$.  Subsequently Crans \cite{Crans2} corrected some
errors in their work and dealt with the cases $n = 2$, $k > 1$.  

\section{Lessons from Homotopy Theory} \label{lessons}

In Grothendieck's famous 600-page letter to Quillen \cite{Gro}, he
proposed developing $n$-category theory as a vast generalization of
homotopy theory, with a special class of $n$-categories --- the
$n$-groupoids --- corresponding to homotopy $n$-types.  When this idea
is finally worked out, we will be able to translate all of homotopy
theory into the language of $n$-groupoids.  Eventually this should
deepen our understanding of the conceptual foundations of homotopy
theory, and help us apply its techniques to other branches of
mathematics.  But even now, with $n$-category theory still in its
squalling infancy, this translation project is worthwhile.  The reason
is that homotopy theory is our best source of insight into
$n$-categories.  We need to be careful here, since homotopy theory
avoids precisely what is most new and interesting about general
$n$-categories, namely the presence of $j$-morphisms that are not
equivalences.  However, this is a bit less of a drawback than it might
at first seem.  After all, the most mysterious aspect of $n$-category
theory is the origin of coherence laws, and these, being implemented as
equivalences, appear already in the context of $n$-groupoids.

\vskip 2em
\begin{center}
{\small
\begin{tabular}{|c|c|}                                    \hline
$\omega$-groupoids &  homotopy types               \\     \hline
$n$-groupoids      &  homotopy $n$-types           \\     \hline
$k$-tuply groupal         &   homotopy types       \\
$\omega$-groupoids        &   of $k$-fold loop spaces \\     \hline
$k$-tuply groupal         &   homotopy $n$-types      \\
$n$-groupoids             &   of $k$-fold loop spaces \\   \hline
$k$-tuply monoidal        &   homotopy types       \\
$\omega$-groupoids        &   of $E_k$ spaces      \\     \hline
$k$-tuply monoidal        &   homotopy $n$-types   \\
$n$-groupoids             &   of $E_k$ spaces      \\     \hline
stable $\omega$-groupoids &   homotopy types       \\
                          &   of infinite loop spaces \\     \hline
stable $n$-groupoids      &   homotopy $n$-types   \\
                          &   of infinite loop spaces \\     \hline 
$\Z$-groupoids            &   homotopy types       \\
                          &   of spectra           \\     \hline
\end{tabular}} \vskip 1em
3.  Translating between $n$-groupoid theory and homotopy theory
\end{center}
\vskip 0.5em

We can begin this translation project with the help of Table 3.  In
explaining this table, we shall use the simplicial approach to
$n$-groupoids instead of the globular approach.  The reason is that
simplicial methods are quite popular among topologists, so all the
necessary machinery has already been developed \cite{May}.  For other
approaches to $n$-groupoids, see the work of Brown and his collaborators
\cite{Brown}, Tamsamani \cite{Tamsamani2}, and Batanin \cite{Batanin}.  

In topology one usually speaks of `$j$-cells' rather than $j$-morphisms.
A `simplicial set' has a set of $j$-cells for each $j \ge 0$, which we
visualize as being shaped like $j$-simplices.  For each $j$-cell $f$
there are certain $(j-1)$-cells $d_0f, \dots, d_nf$ called `faces' and
$(j+1)$-cells $i_0f, \dots, i_{n+1}f$ called `degeneracies'.   One also
requires that the face and degeneracy maps satisfy certain well-known 
relations; we will not need the formulas here. 

A `$j$-dimensional horn' in a simplicial set is, roughly speaking, a
configuration in which all but one of the faces of a $j$-simplex have
been filled in by $(j-1)$-cells in a consistent way.  A simplicial set
for which any horn can be extended to a $j$-cell is called a `Kan
complex'.   A Kan complex is the simplicial version of an
$\omega$-groupoid: a structure like an $n$-groupoid, but without any
cutoff on the dimension of the $j$-morphisms.  

To see how this idea works, suppose we have a Kan complex containing
a `composable' pair of 1-cells $f$ and $g$, meaning that $d_1 f = d_0 g$.
This gives a 2-dimensional horn with $f$ and $g$ as two of its faces, so
we can extend this horn to a 2-cell $F$, which has as its third face
some 1-cell $h$:

\medskip
\begin{center}
\setlength{\unitlength}{0.0100in}%

\begingroup\makeatletter\ifx\SetFigFont\undefined
\def\x#1#2#3#4#5#6#7\relax{\def\x{#1#2#3#4#5#6}}%
\expandafter\x\fmtname xxxxxx\relax \def\y{splain}%
\ifx\x\y   
\gdef\SetFigFont#1#2#3{%
  \ifnum #1<17\tiny\else \ifnum #1<20\small\else
  \ifnum #1<24\normalsize\else \ifnum #1<29\large\else
  \ifnum #1<34\Large\else \ifnum #1<41\LARGE\else
     \huge\fi\fi\fi\fi\fi\fi
  \csname #3\endcsname}%
\else
\gdef\SetFigFont#1#2#3{\begingroup
  \count@#1\relax \ifnum 25<\count@\count@25\fi
  \def\x{\endgroup\@setsize\SetFigFont{#2pt}}%
  \expandafter\x
    \csname \romannumeral\the\count@ pt\expandafter\endcsname
    \csname @\romannumeral\the\count@ pt\endcsname
  \csname #3\endcsname}%
\fi
\fi\endgroup
\begin{picture}(141,120)(228,422)
\thicklines

\put(241,441){\circle*{6}}
\put(360,440){\circle*{6}}
\put(240,440){\line( 3, 4){ 60}}
\put(300,520){\line( 3,-4){ 60}}
\multiput(260,475)(0.50000,0.25000){21}{\makebox(0.4444,0.6667){\SetFigFont{7}{8.4}{rm}.}}
\put(270,480){\line( 0,-1){ 10}}
\put(336,482){\line( 0,-1){ 10}}
\multiput(336,472)(-0.50000,0.25000){21}{\makebox(0.4444,0.6667){\SetFigFont{7}{8.4}{rm}.}}
\multiput(240,440)(11.42857,0.00000){11}{\line( 1, 0){  5.714}}
\multiput(298,445)(0.50000,-0.25000){21}{\makebox(0.4444,0.6667){\SetFigFont{7}{8.4}{rm}.}}
\multiput(308,440)(-0.50000,-0.25000){21}{\makebox(0.4444,0.6667){\SetFigFont{7}{8.4}{rm}.}}
\multiput(300,460)(-0.25000,0.50000){21}{\makebox(0.4444,0.6667){\SetFigFont{7}{8.4}{rm}.}}
\put(300,520){\circle*{6}}
\multiput(300,460)(0.25000,0.50000){21}{\makebox(0.4444,0.6667){\SetFigFont{7}{8.4}{rm}.}}
\multiput(302,479)(0.00000,-10.00000){2}{\line( 0,-1){  5.000}}
\multiput(298,479)(0.00000,-10.00000){2}{\line( 0,-1){  5.000}}

\put(280,465){\makebox(0,0)[lb]{\raisebox{0pt}[0pt][0pt]{$F$}}}
\put(228,431){\makebox(0,0)[lb]{\raisebox{0pt}[0pt][0pt]{$x$}}}
\put(299,527){\makebox(0,0)[lb]{\raisebox{0pt}[0pt][0pt]{$y$}}}
\put(369,433){\makebox(0,0)[lb]{\raisebox{0pt}[0pt][0pt]{$z$}}}
\put(261,481){\makebox(0,0)[lb]{\raisebox{0pt}[0pt][0pt]{$f$}}}
\put(338,484){\makebox(0,0)[lb]{\raisebox{0pt}[0pt][0pt]{$g$}}}
\put(301,422){\makebox(0,0)[lb]{\raisebox{0pt}[0pt][0pt]{$h$}}}
\end{picture}
\end{center}
\medskip

\noindent In this situation, we call $h$ `a composite' of $f$ and $g$,
and think of $F$ as the `process of composition'.  Note that there is
not a unique preferred composite, so composition is not an operation in
the traditional sense.   However, any two composites can be seen to be
equivalent, where two $j$-cells with all the same faces are said to be
`equivalent' if there is a $(j+1)$-cell having them as two of its faces,
the rest being degenerate.    From an algebraic viewpoint, the reason is
that we have defined composition by a {\it universal property}.   

Thanks to the magic of universal properties, Kan complexes are a
wonderfully efficient formalism for studying $\omega$-groupoids. In
particular, there is no need to explicitly list coherence laws!  They
are all implicit in the fact that every horn can be extended to a cell,
and they automatically become explicit if we arbitrarily choose
processes of composition.  For example, given a composable triple of
1-cells, one obtains the associator by considering a 3-simplex with
these as three of its edges and making clever use of the horn-filling
condition.  Likewise, given a composable quadruple of 1-cells, one
obtains the pentagon coherence law by considering a suitable 4-simplex. 
In fact, all the higher coherence laws for associativity, which Stasheff
\cite{Stasheff} organized into polyhedra called `associahedra', have
been obtained from higher-dimensional simplices by Street \cite{Street}
in his simplicial approach to $\omega$-categories.

If we take the liberty of calling Kan complexes `$\omega$-groupoids', we
can set up a correspondence between $\omega$-groupoids and homotopy
types as follows.  Given a topological space $X$, we can form an
$\omega$-groupoid $\Pi(X)$ whose $j$-cells are all the continuous maps
from the standard $j$-simplex into $X$, with faces and degeneracies
defined in the obvious way.  We think of this as the `fundamental
$\omega$-groupoid of $X$'.   Conversely, given an $\omega$-groupoid $G$,
we can form a topological space by taking one geometrical $j$-simplex
for each $j$-cell of $G$ and gluing them all together using the face and
degeneracy maps in the obvious way.  This is called the `geometric
realization of $G$' and denoted $|G|$.    

We thus obtain functors going both ways between the category
$\omega\Gpd$, having Kan complexes as objects and simplicial maps
between these as morphisms, and the category $\Top$, having nice
spaces as objects and continuous maps as morphisms.  (We say a space
is `nice' if it is homotopy equivalent to a CW complex; we use this
nonstandard definition of $\Top$ to exclude pathological spaces like the
topologist's sine curve.)  While these functors are adjoint to one
another, they do not set up an equivalence of categories.  Nonetheless,
we expect that $\omega\Gpd$ and $\Top$ are `the same' in a subtler sense
--- namely, as $\omega$-categories.  More precisely, these categories
should extend to $\omega$-categories, where the 2-morphisms correspond
to homotopies between maps, the 3-morphisms correspond to homotopies
between homotopies, and so on.  The functors
\[ \Pi \maps \Top \to \omega\Gpd, \] 
\[ |\cdot| \maps \omega\Gpd \to \Top\] 
should then extend to $\omega$-functors, giving an equivalence of
$\omega$-categories.

In the absence of this sort of result, topologists have traditionally
used other language to express the fact that $\omega\Gpd$ and $\Top$ are
`the same' for the purposes of homotopy theory.  For example, $\Pi$ and
$|\cdot|$ establish an equivalence between the homotopy category of
$\omega$-groupoids and the homotopy category of nice spaces.  Here the
`homotopy category' is formed by adjoining formal inverses to all maps
inducing isomorphisms of homotopy groups \cite{GZ}, where one defines
the homotopy groups of a Kan complex to be those of its geometric
realization.  An object in the homotopy category of $\Top$ is called a
`homotopy type' \cite{Baues}.

Starting from the correspondence between $\omega$-groupoids and homotopy
types we can set up the other correspondences listed in Table 3.  For
example, we can define an `$n$-groupoid' to be a Kan complex such that
for $j > n + 1$ any configuration in which all the faces of a
$j$-simplex have been filled in by $(j-1)$-cells in a consistent way can
be uniquely extended to a $j$-cell.  This ensures that all cells of
dimension higher than $n$ play the role of equations.   The geometric
realization of an $n$-groupoid is a space with vanishing homotopy groups
above dimension $n$, and the homotopy category of such spaces is called
the category of `homotopy $n$-types'.   The homotopy category of
$n$-groupoids is equivalent to the category of homotopy $n$-types, and
in fact one expects an equivalence of $(n+1)$-categories. 

We define a `$k$-tuply groupal $\omega$-groupoid' to be a Kan complex
with only one $j$-cell for $j < k$.   Under the correspondence between
$\omega$-groupoids and homotopy types, these correspond to homotopy
types with vanishing homotopy groups below dimension $k$.  Similarly, we
define a `$k$-tuply groupal $n$-groupoid' to be an $(n+k)$-groupoid with
only one $j$-cell for $j < k$.  These correspond to homotopy
$(n+k)$-types with vanishing homotopy groups below dimension $k$.  

Recall from the previous section that we can think of an $n$-category
$C$ with only one $j$-morphism for $j < k$  as an $n$-category with
extra structure and properties.  In the globular approach we do this by
`reindexing', constructing an $n$-category whose $j$-morphisms are the
$(j+k)$-morphisms of $C$.   Similarly, it is often useful to think of a
$k$-tuply groupal $n$-groupoid $G$ as an $n$-groupoid with extra
structure and properties.   However, in the simplicial approach we
cannot simply reindex the cells of $G$.   Instead, we can use a standard
construction \cite{May} to form a Kan complex $\Omega^k G$ whose
geometric realization is homotopy equivalent to the $k$th loop space of
the geometric realization of $G$.    This is why in Table 3 we say that
$k$-tuply groupal $n$-groupoids correspond to homotopy $n$-types of
$k$-fold loop spaces.   Similarly, $k$-tuply groupal $\omega$-groupoids
correspond to homotopy types of $k$-fold loop spaces.  

Exactly what extra structure and properties does an $\omega$-groupoid
have if it is $k$-tuply groupal?  In other words, what extra structure
and properties does a space have if it is homotopy equivalent to a
$k$-fold loop space?  This question has inspired the development of many
interesting mathematical tools.  Here we shall describe just one of
these, the `little $k$-cubes operad', invented by Boardman and Vogt
\cite{BV,BV2} and cast into the language of operads by May \cite{May2}.

We begin with the definition of an `operad'.  For each $\ell \ge 0$, an
operad $O$ has a set $O_\ell$ of `$\ell$-ary operations'.  We visualize
such an operation as a tree with one vertex or `node', $\ell$ edges
representing inputs coming in from above, and one edge representing the
output coming out from below:

\vbox{

\begin{center}
\setlength{\unitlength}{0.000500in}%
\begingroup\makeatletter\ifx\SetFigFont\undefined
\def\x#1#2#3#4#5#6#7\relax{\def\x{#1#2#3#4#5#6}}%
\expandafter\x\fmtname xxxxxx\relax \def\y{splain}%
\ifx\x\y   
\gdef\SetFigFont#1#2#3{%
  \ifnum #1<17\tiny\else \ifnum #1<20\small\else
  \ifnum #1<24\normalsize\else \ifnum #1<29\large\else
  \ifnum #1<34\Large\else \ifnum #1<41\LARGE\else
     \huge\fi\fi\fi\fi\fi\fi
  \csname #3\endcsname}%
\else
\gdef\SetFigFont#1#2#3{\begingroup
  \count@#1\relax \ifnum 25<\count@\count@25\fi
  \def\x{\endgroup\@setsize\SetFigFont{#2pt}}%
  \expandafter\x
    \csname \romannumeral\the\count@ pt\expandafter\endcsname
    \csname @\romannumeral\the\count@ pt\endcsname
  \csname #3\endcsname}%
\fi
\fi\endgroup
\begin{picture}(1844,2144)(3279,-3083)
\thicklines
\put(4201,-1861){\circle*{150}}
\put(4201,-3061){\line( 0, 1){1200}}
\put(4201,-1861){\line(-1, 3){300}}
\put(4201,-1861){\line( 1, 3){300}}
\put(4201,-1861){\line( 1, 1){900}}
\put(4201,-1861){\line(-1, 1){900}}
\end{picture}
\end{center}

\medskip
}

\noindent We can compose these trees by feeding the outputs of
$\ell$ of them into one with $\ell$ inputs:

\vbox{

\begin{center}
\setlength{\unitlength}{0.0002500in}%
\begingroup\makeatletter\ifx\SetFigFont\undefined
\def\x#1#2#3#4#5#6#7\relax{\def\x{#1#2#3#4#5#6}}%
\expandafter\x\fmtname xxxxxx\relax \def\y{splain}%
\ifx\x\y   
\gdef\SetFigFont#1#2#3{%
  \ifnum #1<17\tiny\else \ifnum #1<20\small\else
  \ifnum #1<24\normalsize\else \ifnum #1<29\large\else
  \ifnum #1<34\Large\else \ifnum #1<41\LARGE\else
     \huge\fi\fi\fi\fi\fi\fi
  \csname #3\endcsname}%
\else
\gdef\SetFigFont#1#2#3{\begingroup
  \count@#1\relax \ifnum 25<\count@\count@25\fi
  \def\x{\endgroup\@setsize\SetFigFont{#2pt}}%
  \expandafter\x
    \csname \romannumeral\the\count@ pt\expandafter\endcsname
    \csname @\romannumeral\the\count@ pt\endcsname
  \csname #3\endcsname}%
\fi
\fi\endgroup
\begin{picture}(4544,4244)(3279,-5183)
\thicklines
\put(4201,-1861){\circle*{150}}
\put(5701,-1861){\circle*{150}}
\put(7201,-1861){\circle*{150}}
\put(5701,-3661){\circle*{150}}
\put(4201,-3061){\line( 0, 1){1200}}
\put(4201,-1861){\line(-1, 3){300}}
\put(4201,-1861){\line( 1, 3){300}}
\put(4201,-1861){\line( 1, 1){900}}
\put(4201,-1861){\line(-1, 1){900}}
\put(4201,-3061){\line( 5,-2){1500}}
\put(5701,-3661){\line( 0, 1){600}}
\put(5701,-3061){\line( 0, 1){1200}}
\put(5701,-3661){\line( 5, 2){1500}}
\put(7201,-3061){\line( 0, 1){1200}}
\put(7201,-1861){\line(-2, 3){600}}
\put(7201,-1861){\line( 0, 1){900}}
\put(7201,-1861){\line( 2, 3){600}}
\put(5701,-3661){\line( 0,-1){1500}}
\end{picture}
\end{center}

\medskip
}

\noindent In other words, for any $i_1, \dots, i_\ell$ there is a function  
\ban    O_\ell \times O_{i_1} \times \cdots \times
O_{i_\ell} &\to& O_{i_1 + \cdots + i_\ell}  . \\
(f,g_1, \dots, g_\ell) &\mapsto& f \cdot (g_1, \dots , g_\ell)  \ean
We require that composition is `associative', meaning that
\[      f \cdot (g_1 \cdot (h_{11}, \dots, h_{1i_1}), \dots, 
g_\ell \cdot (h_{\ell 1}, \dots, h_{\ell i_\ell})) = \]
\[ (f \cdot (g_1, \dots g_\ell)) \; \cdot \; (h_{11}, \dots, h_{1i_1},
\dots\dots , h_{\ell 1}, \dots, h_{\ell i_\ell}) \]
whenever both sides are well-defined.  
This makes composites such as the following one unambiguous:

\vbox{

\begin{center}

\setlength{\unitlength}{0.000250in}%
\begingroup\makeatletter\ifx\SetFigFont\undefined
\def\x#1#2#3#4#5#6#7\relax{\def\x{#1#2#3#4#5#6}}%
\expandafter\x\fmtname xxxxxx\relax \def\y{splain}%
\ifx\x\y   
\gdef\SetFigFont#1#2#3{%
  \ifnum #1<17\tiny\else \ifnum #1<20\small\else
  \ifnum #1<24\normalsize\else \ifnum #1<29\large\else
  \ifnum #1<34\Large\else \ifnum #1<41\LARGE\else
     \huge\fi\fi\fi\fi\fi\fi
  \csname #3\endcsname}%
\else
\gdef\SetFigFont#1#2#3{\begingroup
  \count@#1\relax \ifnum 25<\count@\count@25\fi
  \def\x{\endgroup\@setsize\SetFigFont{#2pt}}%
  \expandafter\x
    \csname \romannumeral\the\count@ pt\expandafter\endcsname
    \csname @\romannumeral\the\count@ pt\endcsname
  \csname #3\endcsname}%
\fi
\fi\endgroup
\begin{picture}(4980,5444)(3218,-5183)
\thicklines
\put(4201,-1861){\circle*{150}}
\put(5701,-1861){\circle*{150}}
\put(7201,-1861){\circle*{150}}
\put(5701,-3661){\circle*{150}}
\put(3901,-361){\circle*{150}}
\put(4501,-361){\circle*{150}}
\put(5101,-361){\circle*{150}}
\put(6601,-361){\circle*{150}}
\put(7801,-361){\circle*{150}}
\put(3301,-361){\circle*{150}}
\put(7201,-361){\circle*{150}}
\put(4201,-3061){\line( 0, 1){1200}}
\put(4201,-1861){\line(-1, 3){300}}
\put(4201,-1861){\line( 1, 3){300}}
\put(4201,-1861){\line( 1, 1){900}}
\put(4201,-1861){\line(-1, 1){900}}
\put(4201,-3061){\line( 5,-2){1500}}
\put(5701,-3661){\line( 0, 1){600}}
\put(5701,-3061){\line( 0, 1){1200}}
\put(5701,-3661){\line( 5, 2){1500}}
\put(7201,-3061){\line( 0, 1){1200}}
\put(7201,-1861){\line(-2, 3){600}}
\put(7201,-1861){\line( 0, 1){900}}
\put(7201,-1861){\line( 2, 3){600}}
\put(5701,-3661){\line( 0,-1){1500}}
\put(3301,-961){\line( 0, 1){600}}
\put(3901,-961){\line( 0, 1){600}}
\put(4501,-961){\line( 0, 1){600}}
\put(5101,-961){\line( 0, 1){600}}
\put(6601,-961){\line( 0, 1){600}}
\put(7201,-1036){\line( 0, 1){675}}
\put(7801,-961){\line( 0, 1){600}}
\put(5101,-361){\line( 1, 4){150}}
\put(5101,-361){\line(-1, 4){150}}
\put(4501,-361){\line( 0, 1){600}}
\put(3901,-361){\line( 1, 2){300}}
\put(3901,-361){\line( 0, 1){600}}
\put(3901,-361){\line(-1, 2){300}}
\put(6601,-361){\line( 1, 2){300}}
\put(6601,-361){\line( 0, 1){600}}
\put(6601,-436){\line(-2, 5){274.138}}
\put(7801,-361){\line(-3, 5){363.971}}
\put(7801,-361){\line(-1, 4){150}}
\put(7801,-361){\line( 1, 4){150}}
\put(7801,-361){\line( 3, 5){363.971}}
\end{picture}
\end{center}

\medskip
}

\noindent We also require the existence of an unit operation $1 \in O_1$
such that
\[     1 \cdot (f) = f, \qquad f \cdot (1,\dots, 1) = f \]
for all $f \in O_\ell$.  

What we have so far is an {\it planar operad}.   For a
full-fledged operad, we also assume that there are right actions of
the symmetric groups $S_\ell$ on the sets $O_\ell$ for which the
following compatibility conditions hold.  First, for any $f \in O_\ell$,
$\sigma \in S_\ell$, and $g_j \in O_{i_j}$ for $1 \le j \le \ell$, we
require
\[     (f\sigma) \cdot (g_{\sigma(1)}, \dots, g_{\sigma(\ell)} ) = 
(f \cdot (g_1, \dots, g_\ell))\, \rho(\sigma), \]
where
\[      \rho \maps S_\ell \to S_{i_1 + \cdots + i_\ell}  \]
is the obvious homomorphism.  
Second, for any $f \in O_\ell$, and $g_j \in O_{i_j}$, $\sigma_j
\in S_{i_j}$ for $1 \le j \le \ell$, we require
\[     f \cdot (g_1\sigma_1, \dots, g_\ell \sigma_\ell) = 
(f \cdot (g_1, \dots, g_\ell))\, \rho'(\sigma_1,\dots,\sigma_\ell), \]
where
\[      \rho' \maps S_{i_1} \times \cdots \times S_{i_\ell}
\to S_{i_1 + \cdots + i_\ell} \]
is the obvious homomorphism. 

Just as groups are interesting for their actions, operads are 
interesting for their `algebras'.  Given an operad $O$, an
`$O$-algebra' is a set $A$ equipped with actions
\[      \alpha \maps O_\ell \times A^\ell \to A, \] 
or equivalently, maps
\[      \alpha \maps O_\ell \to \hom(A^\ell,A) \] 
representing the $\ell$-ary operations of $O$ as actual operations on
$A$.  We require that $\alpha$ sends the identity operation $1
\in O_1$ to the identity function on $A$ and sends composites to
composites: 
\[      \alpha(f \cdot
(g_1, \dots, g_\ell)) = \alpha(f) \circ (\alpha(g_1) \times \cdots \times
\alpha(g_\ell)) .\] 
We also require that
\[        \alpha(f\sigma) = \alpha(f)\sigma, \] 
where $f \in O_\ell$, and $\sigma \in S_\ell$ acts on $\hom(A^\ell,A)$ on the
right by permuting the factors in $A^\ell$.  We omit this requirement if
$O$ is merely planar.

More generally, one can define operads and their algebras in any
symmetric monoidal category $C$, which amounts to replacing sets and
functions in the above definitions by objects and morphisms in $C$, and
replacing the Cartesian product by the tensor product in $C$.  We shall
mainly be interested in operads in the category $\Top$, which are called
`topological operads'.  Spaces with extra structure and properties can 
often be described as algebras of topological operads.

The most interesting example for us is the `little $k$-cubes operad', 
$C(k)$.   Here the space $C(k)_\ell$ of $\ell$-ary operations is the space 
of $\ell$ disjoint $k$-cubes linearly embedded in the cube $[0,1]^k$ 
with their edges parallel to the coordinate axes:

\vbox{

\begin{center}

\setlength{\unitlength}{0.000500in}%
\begingroup\makeatletter\ifx\SetFigFont\undefined
\def\x#1#2#3#4#5#6#7\relax{\def\x{#1#2#3#4#5#6}}%
\expandafter\x\fmtname xxxxxx\relax \def\y{splain}%
\ifx\x\y   
\gdef\SetFigFont#1#2#3{%
  \ifnum #1<17\tiny\else \ifnum #1<20\small\else
  \ifnum #1<24\normalsize\else \ifnum #1<29\large\else
  \ifnum #1<34\Large\else \ifnum #1<41\LARGE\else
     \huge\fi\fi\fi\fi\fi\fi
  \csname #3\endcsname}%
\else
\gdef\SetFigFont#1#2#3{\begingroup
  \count@#1\relax \ifnum 25<\count@\count@25\fi
  \def\x{\endgroup\@setsize\SetFigFont{#2pt}}%
  \expandafter\x
    \csname \romannumeral\the\count@ pt\expandafter\endcsname
    \csname \romannumeral\the\count@ pt\endcsname
  \csname #3\endcsname}%
\fi
\fi\endgroup
\begin{picture}(3624,3624)(3589,-5173)
\thicklines
\put(3601,-1561){\line( 0,-1){3600}}
\put(3601,-5161){\line( 1, 0){3600}}
\put(3601,-1561){\line( 1, 0){3600}}
\put(7201,-1561){\line( 0,-1){3600}}
\put(4201,-2161){\makebox(6.6667,10.0000){\SetFigFont{10}{12}{rm}.}}
\put(4201,-1861){\line( 0,-1){900}}
\put(4201,-2761){\line( 1, 0){900}}
\put(5101,-2761){\line( 0, 1){900}}
\put(5101,-1861){\line(-1, 0){900}}
\put(6001,-2461){\line( 0,-1){600}}
\put(6001,-3061){\line( 1, 0){600}}
\put(6601,-3061){\line( 0, 1){600}}
\put(6001,-2461){\line( 1, 0){600}}
\put(6001,-3361){\line(-1, 0){1200}}
\put(4801,-3361){\line( 0,-1){1200}}
\put(4801,-4561){\line( 1, 0){1200}}
\put(6001,-4561){\line( 0, 1){1200}}
\put(6901,-4561){\line(-1, 0){300}}
\put(6601,-4561){\line( 0,-1){300}}
\put(6601,-4861){\line( 1, 0){300}}
\put(6901,-4861){\line( 0, 1){300}}
\end{picture}

\end{center}
\medskip
}

An algebra of the little $k$-cubes operad is called an `$E_k$ space'.
Given a space $X$ with a basepoint, there is an obvious way to make 
$\Omega^k X$ into an $E_k$ space using the fact that it consists of maps 
$f \maps [0,1]^k \to X$ sending the boundary of the cube to the basepoint.   
Conversely, the set of components of an $E_k$ space automatically
becomes a monoid, and if this monoid is a group then the $E_k$ space is
homotopy equivalent to a $k$-fold loop space.  Homotopy types of $E_k$
spaces are thus a bit more general than homotopy types of $k$-fold loop
spaces.  While the latter correspond to $k$-tuply groupal
$\omega$-groupoids, the former should correspond to `$k$-tuply monoidal
$\omega$-groupoids', that is, $\omega$-categories with only one
$j$-morphism for $j < k$, for which all $j$-morphisms are equivalences
for $j > k$ --- but not necessarily for $j = k$.    There should
similarly be a correspondence between $k$-tuply monoidal $n$-groupoids
and homotopy $n$-types of $E_k$ spaces. 

At the end of the previous section we listed various processes going
between neighboring entries of Table 2.  When restricted to $k$-tuply
groupal $n$-groupoids, most of these processes have well-known
topological interpretations, which we summarize below.  We let
$n\Type_k$ stand for the category of homotopy $n$-types of $k$-fold loop
spaces.  It is easiest to define this as the full subcategory of the
homotopy category of nice spaces $X$ with basepoint such that
$\pi_j(X) = 0$ unless $k \le j \le n+k$.  However, by repeated
looping, we shall think of the objects of $n\Type_k$ as $k$-fold loop
spaces with vanishing homotopy groups above dimension $n$.  A technical
point worth noting is that for $k = 0$, these really correspond to
$n$-groupoids {\it equipped with a distinguished 0-cell}.

\vskip 1em $\bullet$ \noindent {\it Decategorification: $(n,k) \to
(n-1,k)$}.   Let $X$ be a $k$-fold loop space with $\pi_j(X) = 0$ for $j
> n$. Then we can kill off its $n$th homotopy group by attaching cells,
obtaining a $k$-fold loop space $\Decat X$ with $\pi_j(X) = 0$ for $j >
n - 1$.  There is a map $X \to \Decat X$ inducing isomorphisms
$\pi_j(X) \iso \pi_j(\Decat X)$ for $j \le n - 1$. This process gives a
functor $\Decat \maps n\Type_k \to (n-1)\Type_k$.

\vskip 1em $\bullet$ \noindent
{\it Discrete categorification: $(n,k) \to (n+1,k)$}.  
The forgetful functor $\Disc \maps n\Type_k \to (n+1)\Type_k$ is the left
adjoint of $\Decat$.  

\vskip 1em $\bullet$ \noindent
{\it Delooping: $(n,k) \to (n+1,k-1)$}.   
There is a $(k-1)$-fold loop space $BX$, called the `classifying
space' or `delooping' of $X$, with $\pi_j(BX) \iso \pi_{j-1}(X)$.
This process gives a functor $B \maps n\Type_k \to (n+1)\Type_{k-1}$.

\vskip 1em $\bullet$ \noindent
{\it Looping: $(n,k) \to (n-1,k+1)$}.  
For $k > 0$, there is a $(k+1)$-fold loop space $\Omega X$, consisting
of all based loops in $X$, with $\pi_j(\Omega X) \iso \pi_{j-1}(X)$.
This process gives a functor $\Omega \maps n\Type_k \to
(n-1)\Type_{k+1}$ that is right adjoint to $B$.   For $k > 0$, delooping
followed by looping is naturally isomorphic to the identity functor.

\vskip 1em $\bullet$ \noindent
{\it Forgetting monoidal structure: $(n,k) \to (n,k-1)$}.  
There is a forgetful functor $F \maps n\Type_k \to n\Type_{k-1}$.  

\vskip 1em $\bullet$ \noindent
{\it Stabilization: $(n,k) \to (n,k+1)$}. 
There is a left adjoint to $F$ called `stabilization', which
we denote by $S \maps n\Type_k \to n\Type_{k+1}$.  If we define the
`suspension' functor $\Sigma \maps (n-1)\Type_k \to n\Type_k$ to be the
left adjoint of the composite
\[
\begin{diagram} [(n-1)\Type_k]
\node{n\Type_k} \arrow{e,t}{\Omega} 
\node{(n-1)\Type_{k+1}} \arrow{e,t}{F}
\node{(n-1)\Type_k}
\end{diagram}
\]
then stabilization is naturally isomorphic to suspension followed by
looping.  The Freudenthal suspension theorem says that stabilization is
an equivalence for $k \ge n + 2$. 

\vskip 1em $\bullet$ \noindent
{\it Forming the generalized center: $(n,k) \to (n,k+1)$}.  
This process needs to be adapted to stay within the world of $k$-tuply
groupal groupoids, or in other words, homotopy $n$-types of $k$-fold
loop spaces.  We let $\Aut(X)$ be the automorphism group of $X$ as an
object of $n\Type_k$, and define $ZX = \Omega^k(\Aut(X))$.  
\vskip 1em

The fact that stabilization is an equivalence for $k \ge n + 2$ leads us
to define a `stable $n$-groupoid' to be an $(n+2)$-tuply groupal
$n$-groupoid.   We expect, in fact, that there is an $(n+1)$-category of
stable $n$-groupoids, and that the decategorification of this is the
$n$-category of stable $(n-1)$-groupoids.  If we could take the inverse
limit in a suitable sense, we would hope to obtain an $\omega$-category
of `stable $\omega$-groupoids'.   These should correspond to what
topologists call `infinite loop spaces', an infinite loop space being a
sequence of spaces $X_0,X_1,\dots$, equipped with homeomorphisms
$f_k \maps X_k \to \Omega X_{k+1}$.  

Infinite loop spaces play an important role in stable homotopy theory
\cite{Adams,EKMM,May2}.  A closely related concept is that of a
`spectrum'.  Just as an infinite loop space should correspond to a
stable $\omega$-groupoid, a spectrum should correspond to a
`$\Z$-groupoid', some sort of gadget with $j$-morphisms for all $j \in
\Z$, all of which are equivalences.  In fact, strict $\Z$-categories are
easily defined in the globular approach, and strict $\Z$-groupoids then
work out to be the same as $\Z$-graded chain complexes of abelian
groups.  Spectra can be viewed as a generalization of such chain
complexes.  Apart from what we know about spectra, however, the theory
of `weak $\Z$-categories' remains largely {\it terra incognita}.

\section{Examples of Categorification}  \label{examples}

In what follows, we consider iterated categorifications of the natural
numbers, and some variations on this theme where we adjoin formal
inverses or duals.  We shall see that quite a bit of mathematics amounts
to the study of the resulting objects, some of which are astoundingly
complicated.  The more `$n$-groupoidal' examples are familiar from
homotopy theory, but the importance of some of the more
`$n$-categorical' ones has only become clear in recent work on
topological quantum field theory.

The natural numbers are the free monoid on one element, while the
integers are the free group on one element.  Some categorified analogs
of these notions are listed in Table 4.   In the rest of this section,
we work through these examples in detail.  

We begin with the simplest algebraic structure of all: sets.   The free
set on one element is just the one-element set, denoted by 1.  Note that
here we are using the word `the' in a generalized sense.  Since all sets
with the same cardinality are isomorphic, any singleton $\{x\}$ can be
regarded as `the' 1-element set.    In a set, when we speak of `the'
element with some property, we imply that any other element having this
property is {\it equal} to this one, but in a category, when we speak of
`the' object with some property, we mean that any other object having
this property is {\it isomorphic} to this one --- typically by means of
a uniquely specified isomorphism.  More generally, when we speak of
`the' object of an $n$-category having some property, we mean that any
other object with this property is {\it equivalent} --- typically by
means of an equivalence which is specified uniquely up to an 
equivalence which is specified uniquely up to $\dots$ and so on.  This
recursive weakening of the notion of uniqueness, and therefore of the
meaning of `the', is fundamental to categorification.  

\vskip 2em
\begin{center}
{\small
\begin{tabular}{|c|c|}                    \hline
sets     &  1                                   \\     \hline
monoids  &  $\N$                                \\     \hline
groups   &  $\Z$                                \\     \hline
$k$-tuply monoidal        &   $n\Braid_k$       \\
$n$-categories            &                     \\     \hline
$k$-tuply monoidal        &   $\Braid_k$        \\
$\omega$-categories       &                     \\     \hline
stable $n$-categories     &   $n\Braid$         \\     \hline
stable $\omega$-categories&   $\Braid_\infty$   \\     \hline
$k$-tuply monoidal        &   $n\Tang_k$        \\ 
$n$-categories with duals &                     \\     \hline
stable $n$-categories     &   $n\Cob$           \\      
with duals                &                     \\     \hline  
$k$-tuply groupal         &   $\Pi_n(\Omega^k S^k)$  \\
$n$-groupoids             &                          \\     \hline
$k$-tuply groupal         &   $\Pi(\Omega^k S^k)$    \\
$\omega$-groupoids        &                          \\     \hline
stable $\omega$-groupoids &   $\Pi(\Omega^\infty S^\infty)$ \\   \hline
$\Z$-groupoids            &   the sphere spectrum           \\   \hline
\end{tabular}} \vskip 1em
4.  Algebraic structures and the free such structures on one generator
\end{center}
\vskip 0.5em

Next consider the free category on one object.  This is just the 
category with one object $x$ and one morphism $1_x \maps x \to x$.  More 
generally, for each $n$, consider the free $n$-category on one
object.  In the globular approach to $n$-categories we may take this to
be the $n$-category with one object $x$, one morphism $1_x \maps x \to
x$, one 2-morphism $1_{1_x} \maps 1_x \To 1_x$, and so on up to one
$n$-morphism.  This is an $n$-groupoid, namely the globular version
of the fundamental $n$-groupoid of a point.   By the remarks in the
previous paragraph, the fundamental $n$-groupoid of any contractible
space may also be considered `the free $n$-category on one object'. 

Things become more interesting when we generalize further and consider
the free $k$-tuply monoidal $n$-category on one object, which we denote
by $n\Braid_k$, for reasons soon to be apparent.   Let us see what this
looks like in the simplest cases, namely $n = 0$ and $1$.  The
Stabilization Hypothesis says that we only need to consider $k \le n +
2$.

\vskip 1em $\bullet$ \noindent
$0\Braid_0$, the free set on one element.  This is the one-element
set, 1.

\vskip 1em $\bullet$ \noindent
$0\Braid_1$, the free monoid on one element.  This is the
natural numbers, $\N$, with addition as its monoid structure.

\vskip 1em $\bullet$ \noindent
$0\Braid_2$, the free commutative monoid on one element.  This
is again $\N$, now regarded as a commutative monoid with addition as
its monoid structure. 

\vskip 1em $\bullet$ \noindent
$1\Braid_0$, the free category on one object.  This is the
category with one object and one morphism.

\vskip 1em $\bullet$ \noindent
$1\Braid_1$, the free monoidal category on one object $x$.  The objects
of this category are the tensor powers $x^{\tensor \ell}$, and the only
morphisms are identity morphisms.  (Here we are using the Mac Lane
coherence theorem to make $1\Braid_1$ into a strict monoidal category.)
This is the discrete categorification of $0\Braid_1$.

\vskip 1em $\bullet$ \noindent
$1\Braid_2$, the free braided monoidal category on one object $x$.  This
is the braid groupoid.  The objects of this groupoid are the tensor
powers $x^{\tensor \ell}$, and the only morphisms are automorphisms,
with $\End(x^{\tensor \ell})$ being the $\ell$-strand braid group
$B_\ell$, which has generators $\sigma_i$ ($1 \le i \le \ell - 1$) and
relations
\ban        
\sigma_i \sigma_{i+1} \sigma_i &=& \sigma_{i+1} \sigma_i \sigma_{i+1}, \\
\sigma_i \sigma_j &=& \sigma_j \sigma_i  \qquad\qquad\qquad
{\rm for} \; |i - j| > 1. \ean

\vskip 1em $\bullet$ \noindent
$1\Braid_3$, the free symmetric monoidal category on one object $x$.
This is the symmetric groupoid.  The objects of this groupoid are the
tensor powers $x^{\tensor \ell}$, and the only morphisms are automorphisms,
with $\End(x^{\tensor \ell})$ being the symmetric group on $\ell$
letters, $S_\ell$.   This group can be given a presentation like that
of $B_\ell$, but with the additional relations
\[    \sigma_i^2 = 1 .\]
The symmetric groupoid is equivalent to the category with finite sets as
objects and bijections as morphisms.   Note that, like 
the category $\FinSet$ having finite sets as objects and arbitrary
functions as morphisms, this category is a rig category with $\N$ as 
its decategorification.
\vskip 1em

The following heuristic argument allows us to guess a general formula
for $n\Braid_k$.  In all the cases considered above, $n\Braid_k$ is a
$k$-tuply monoidal $n$-groupoid.  This should hold in general, since for
$j > 0$ the $j$-morphisms of $n\Braid_k$ arise from coherence laws, and
should thus be equivalences.  We also expect that $n\Braid_k$ is the
`free $k$-tuply monoidal $n$-groupoid on one object'.  By Table 3,
$n\Braid_k$ should thus be the fundamental $n$-groupoid of some $E_k$
space $X_{n,k}$ with vanishing homotopy groups above dimension $n$.  We
also expect that $\Decat((n+1)\Braid_k) \simeq n\Braid_k$, so there
should be some sort of inverse limit, the `free $k$-tuply monoidal
$\omega$-groupoid on one object', which we denote by $\Braid_k$.
Corresponding to this there should be a sequence of maps $X_{n+1,k} \to
X_{n,k}$, with $X_{n,k}$ obtained from $X_{n+1,k}$ by killing its
$(n+1)$st homotopy group, and with the inverse limit being a space $X_k$
whose fundamental $\omega$-groupoid is $\Braid_k$.  By the
correspondence between $k$-tuply monoidal $\omega$-groupoids and $E_k$
spaces, we expect that $X_k$ is the `free $E_k$ space on one point'.

While this argument involves many forms of reasoning that have not yet
been made rigorous, there is at least a precise meaning to the `free
$E_k$ space on one point'.  Given any operad $O$ in the category
of pointed spaces and any pointed space $X$, there is a standard
construction of the `free $O$-algebra on $X$', due to May \cite{May2}.  
To form the free $O$-algebra on one point, we should first form the free
free pointed space on one point, namely $S^0$, and then apply
this standard construction.  The result is the disjoint union 
\[  \coprod_{\ell = 0}^\infty O_\ell /S_\ell \]
which becomes an $O$-algebra in a tautologous way.  

In the case at hand, since the the operad for $E_k$ spaces is the little
$k$-cubes operad $C(k)$, the free $E_k$ space on one point is
\[             X_k = \coprod_{\ell = 0}^\infty C(k)_\ell/S_\ell . \]
We expect, therefore, an equivalence
\[ n\Braid_k \simeq \Pi_n(X_k) . \]
To work with this equivalence, it is helpful to note that $C(k)_\ell$ is 
homotopy equivalent to the `configuration space' \cite{Segal} of $\ell$ 
distinct points in the $k$-cube:
\[      \{(x_1,\dots,x_\ell) \in [0,1]^k \;\colon \; x_i \ne x_j
{\rm \; if\;} i \ne j \}.      \] 
Moreover, this homotopy equivalence is compatible with the obvious
actions of $S_\ell$.  It follows that $X_k$ is homotopy equivalent to
the pointed space of all finite sets of distinct points in the $k$-cube,
where the empty set plays the role of basepoint.

To see how this works in an example, consider the case $n = 1$, $k = 2$.
The space $X_2$ is equivalent to the space of finite sets of points in
the square.  A path in $X_2$ amounts to a braid with an arbitrary number
of strands.  It follows that the fundamental groupoid of $X_2$ is
equivalent to the braid groupoid, so
\[     1\Braid_2 \simeq   \Pi_1(X_2)  .\]
More generally, in the globular approach to $n$-categories, the
$n$-morphisms of $n\Braid_k$ should correspond to certain
$n$-dimensional surfaces in $[0,1]^{n+k}$, which we could call
`$n$-braids in $n+k$ dimensions'.  We see here an instance of a general
theme, namely that in a $k$-tuply monoidal $n$-category the number $n$
often plays the role of `dimension', while $k$ plays the role of
`codimension'.

Now let us turn to the `free $k$-tuply groupal $n$-groupoid on one
object', which we temporarily denote by $G_{n,k}$.   For $n = 0$, $k =
1$ this is just the free group on one element, namely the integers,
$\Z$.   For higher values of $n$ and $k$ we may thus regard $G_{n,k}$ as
a categorified, stabilized version of $\Z$.   In what follows we restrict
attention to the case $k > 0$, since in this case a $k$-tuply groupal
$n$-groupoid automatically has a distinguished object, the unit.  

Again we can use a heuristic argument to guess a formula for $G_{n,k}$. 
Since we expect that $G_{n,k} \simeq \Decat G_{n+1,k}$, there should be
some sort of inverse limit $G_k$, the `free $k$-tuply groupal
$\omega$-groupoid on one object'.  By Table 3 we expect $G_k$ to be the
fundamental $\omega$-groupoid of some $k$-fold loop space, so the
problem is to determine this space.  Just as the group $\Z$ is obtained 
from the monoid $\N$ by adjoining formal inverses, $G_k$ should be 
obtained from $\Braid_k$ by adjoining formal weak inverses for all 
objects.  More generally, we should be able to turn any $k$-tuply
monoidal $\omega$-groupoid into a $k$-tuply groupal $\omega$-groupoid by
adjoining formal weak inverses of objects.   In the language of homotopy
theory, this process should turn $E_k$ spaces into $k$-fold loop spaces.
In fact, this process is familiar in homotopy theory under the name of
`group completion' \cite{BP}.   Since $n\Braid_k$ is the fundamental
$\omega$-groupoid of $X_k$, we thus expect $G_k$ to be the fundamental
$\omega$-groupoid of the group completion of $X_k$.  The group
completion of $X_k$ is homotopy equivalent to $\Omega^k S^k$, so we
expect an equivalence 
\[           G_{n,k} \simeq \Pi_n(\Omega^k S^k)  .\] 
Group completion automatically gives a map from $X_k$ to
$\Omega^k S^k$ which induces a $k$-tuply monoidal $n$-functor 
\[           \Pi_n(X_k) \to \Pi_n(\Omega^k S^k)   .\]

If the above guesses are correct, $n$-braids and the homotopy types of
spheres play a fundamental role in $n$-category theory.  Since the 
homotopy groups of spheres are notoriously hard to compute, this means
that $n$-category theory has a certain built-in complexity.  Perhaps we
should amplify on this a bit.  Suppose that $C$ is a globular
$(n+k)$-category and $x$ is any object of $C$.  Let $1_1 = 1_x$ and
recursively define $1_{i+1}$ to be $1_{1_i}$.   Then we are claiming
that any $k$-morphism $f \maps 1_{k-1} \to 1_{k-1}$ determines, at least
up to equivalence, a $k$-tuply monoidal $n$-functor
\[           n\Braid_k \to \End(1_{k-1})     \]
mapping the generator of $n\Braid_k$ to $f$.  Moreover, if $f$ is 
an equivalence, we claim this factors through a $k$-tuply groupal
$n$-functor
\[          \Pi_n(\Omega^k S^k) \to \End(1_{k-1})  \]
whose range consists entirely of equivalences.  

In our previous exploration of these ideas \cite{BD}, we emphasized the
importance of a notion lying halfway between $n$-groupoids and fully
general $n$-categories, which we called `$n$-categories with duals'.
The idea here is that duals are an interesting generalization of
inverses.  In particular, the Tangle Hypothesis states that there is a
`free $k$-tuply monoidal $n$-category with duals on one object',
$n\Tang_k$, having as $n$-morphisms certain $n$-dimensional surfaces in
$[0,1]^{n+k}$ called `framed $n$-tangles in $n+k$ dimensions'.  For $n =
1, k = 2$ this was proved by Freyd and Yetter, Turaev, and Shum
\cite{FY,Shum,T,Y}.  This special case serves as the basis of recent
work on 3-dimensional topological quantum field theory.  Indeed, the
fact that $n$-categories with duals are more general than $n$-groupoids
is the reason why topological quantum field theory can give more refined
information than homotopy theory.  Recently progress has been made on
the case $n = 2$, $k = 2$, which has also illuminated the theory of
2-braids in 4 dimensions \cite{BL,BL2,CRS,KT}.  In general, we expect
that in the stable range $n\Tang_k$ is equivalent to the stable
$n$-category of `framed cobordisms'.  Also, the universal property of
$n\Tang_k$ should give a $k$-tuply monoidal $n$-functor
\[ T \maps n\Tang_k \to \Pi_n(\Omega^k S^k), \] 
generalizing the Thom-Pontryagin construction.  For more details the
reader must turn to the references.  Our main point here is that if the
Tangle Hypothesis holds, a great deal, not only of homotopy theory, but
also of topological quantum field theory arises naturally from
the study of categorified analogs of $\Z$!

To conclude, let us note that all the entries in Table 4 should be
equipped with `multiplication' as well as `addition' operations.  For
example, by virtue of being the free monoid on one element, $\N$
automatically becomes a rig in the following way: given an element $n
\in \N$, there is a unique monoid homomorphism $f \maps \N \to \N$ with
$f(1) = n$, namely multiplication by $n$.  Likewise, by virtue of being
the free group on one element, $\Z$ automatically becomes a ring.
Categorifying once and stabilizing various numbers of times, we see that
for $k > 0$, $1\Braid_k$ is a rig category and $\Pi_1(\Omega^k S^k)$ is
a `ring category': a rig category for which objects have additive
inverses.  This pattern should continue throughout the rest of Table 4.

For example, in stable homotopy theory it is well known that the sphere
spectrum is a `ring spectrum' \cite{MQRT}.  We expect that the sphere
spectrum corresponds to the `free $\Z$-groupoid on one object'.  Indeed,
Joyal has called the sphere spectrum `the true integers', since it is an
infinitely categorified, infinitely stabilized analog of $\Z$.

\section{$n$-Categories and the Algebra of Opetopes}   \label{opetopes}

We hope the previous sections have whetted the reader's appetite for a
rigorous theory of $n$-categories by sketching a bit of what we might do
if we had one.  Now we turn to the issue of actually developing this
theory.  It seems that any definition of $n$-category involves a choice
of the basic shapes of $j$-morphisms --- globes, simplices, or whatever.
It also involves a choice of ways to compose $j$-morphisms by gluing
these basic shapes together.  Most importantly, it requires a careful
treatment of coherence laws.  In what follows we present an approach in
which all these issues are handled simultaneously using the formalism of
operads.  In this approach, the basic shapes of $j$-morphisms are the
$j$-dimensional `opetopes'.  The allowed ways of composing $j$-morphisms
correspond precisely to the $(j+1)$-dimensional opetopes.  Moreover, the
sequence of higher coherence laws satisfied by composition correspond to
opetopes of ever higher dimension.

Before going into the details, let us sketch how this works.  First
consider some low-dimensional opetopes.  The only 0-dimensional opetope
is the point:

\medskip
\centerline{\epsfysize=0.08in\epsfbox{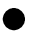}}
\medskip

\noindent The only way to glue together 0-dimensional opetopes is the
trivial way: the identity operation.  The only 1-dimensional opetope is
thus the interval, or more precisely the arrow:

\medskip
\centerline{\epsfysize=0.12in\epsfbox{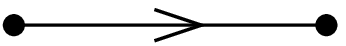}}
\medskip

\noindent The allowed ways of gluing together 1-dimensional opetopes are
given by the 2-dimensional opetopes.  The first few 2-dimensional
opetopes are as follows:

\medskip
\centerline{\epsfysize=0.8in\epsfbox{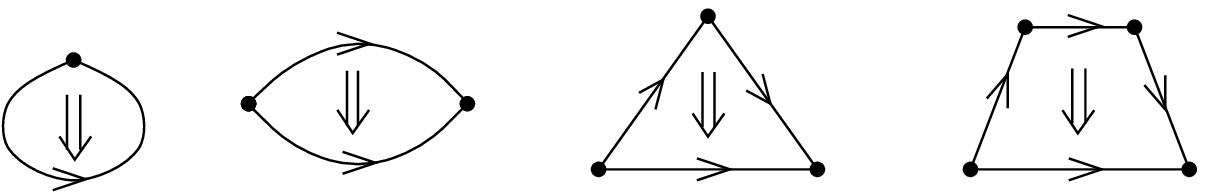}}
\medskip

\noindent  The allowed ways of gluing together 2-dimensional opetopes are
given by the 3-dimensional opetopes.  There are many of these; a simple
example is as follows:

\medskip
\centerline{\epsfysize=1.2in\epsfbox{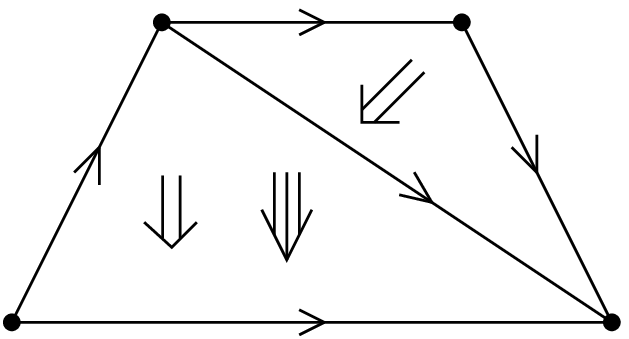}}
\medskip

\noindent This may be a bit hard to visualize, but it depicts a
3-dimensional shape whose front consists of two 3-sided `infaces', and
whose back consists of a single 4-sided `outface'.  We have drawn double
arrows on the infaces but not on the outface.  Note that while this
shape is topologically a ball, it cannot be realized as a polyhedron
with planar faces.  This is typical of opetopes.   

In general, an $(n+1)$-dimensional opetope has any number of infaces and
exactly one outface: the infaces are $n$-dimensional opetopes glued
together in a tree-like pattern, while the outface is a single
$n$-dimensional opetope.  For example, the 3-dimensional opetope above
corresponds to the following tree:

\medskip
\centerline{\epsfysize=1.2in\epsfbox{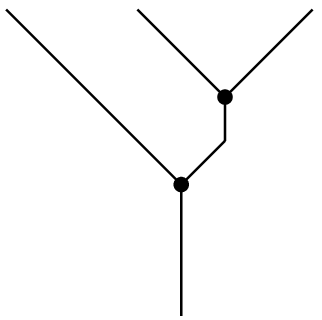}}
\medskip

\noindent The two triangular infaces of the opetope correspond to the
two nodes in this tree.  This is a rather special tree; in general, we
allow nonplanar trees with any number of nodes and any number of edges
coming into each node.

Our approach to $n$-categories is a bit like the Kan complex approach to
$n$-groupoids described in Section \ref{lessons}, but with simplicial
sets replaced by `opetopic sets'.  Basically, an opetopic set is a set
of `cells' shaped like opetopes, such that any face of a cell is again
a cell.  In an $n$-category, the $j$-dimensional cells play the role of
$j$-morphisms.  An opetopic set is an $n$-category if it satisfies
the following two properties:

1) {\it ``Any niche has a universal occupant.''}  A `niche' is
a configuration where all the infaces of an opetope have been filled in
by cells, but not the outface or the opetope itself:

\medskip
\centerline{\epsfysize=1.5in\epsfbox{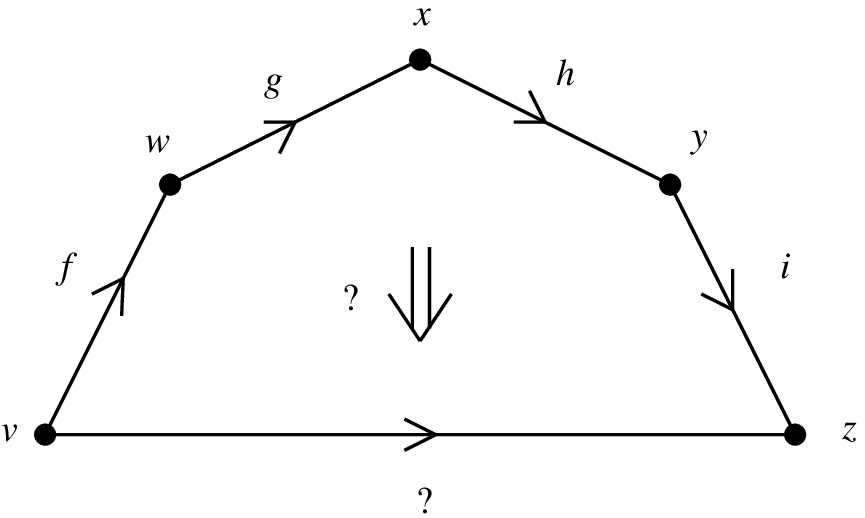}}
\medskip

\noindent An `occupant' of the niche is a way of extending this
configuration by filling in the opetope (and thus its ouface) with a
cell:

\medskip
\centerline{\epsfysize=1.5in\epsfbox{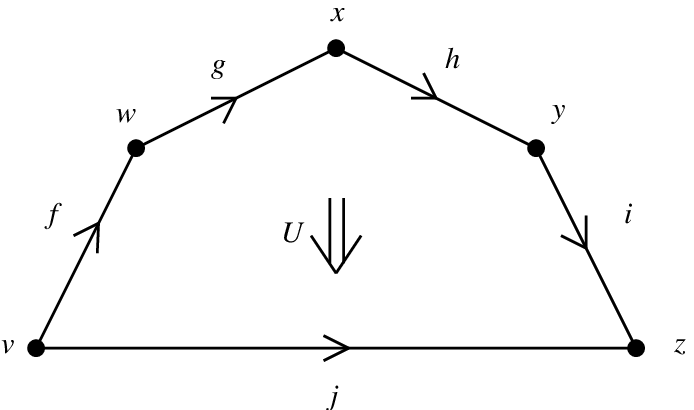}}
\medskip

\noindent The `universality' of an occupant means roughly that every other
occupant factors through the given one {\it up to equivalence}.  To make
this precise we need to define universality in a rather subtle recursive
way.  We may think of a universal occupant of a niche as `a process of
composing' the infaces, and its outface as `a composite' of the infaces.

2) {\it ``Composites of universal cells are universal.''}  Suppose
that $U,V,$ and $W$ below are universal cells:

\medskip
\centerline{\epsfysize=1.5in\epsfbox{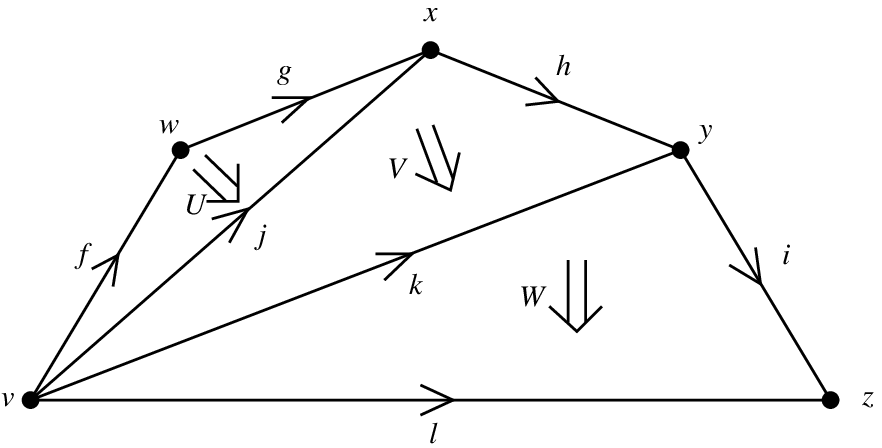}}
\medskip

\noindent Then we can compose them, and we are guaranteed that their
composite is again universal, and thus that the outface $l$ is a composite
of the cells $f,g,h,i$.  Note that a process of composing $U,V,W$
is described by a universal occupant of a niche of one higher
dimension.  

Note that in this approach to $n$-categories, composition of cells is
not an operation in the traditional sense: the composite is defined by a
universal property, and is thus unique only up to equivalence.   Only at
the top level, for the $n$-cells of an $n$-category, is the composite
truly unique.  The main advantage of defining composition by a universal
property is that we do not need to list coherence laws: all the right
coherence laws arise automatically!    This is a very important point,
because in some sense it answers the puzzle concerning the origin of
coherence laws. 

At first this answer may seem as puzzling as the puzzle it answers.
{\it Why} does defining composition by a universal property
automatically generate all the right coherence laws?  One reason is that
coherence laws are `right' when they hold in interesting examples, and
in these examples composition is usually defined by a universal
property.  Consider for example the categorified version of $\N$
discussed in the Introduction: the category $\FinSet$.  Corresponding to
addition in $\N$, the category $\FinSet$ has finite coproducts, i.e.,
disjoint unions.  Coproducts are defined by a universal property, and
this universal property immediately implies number of things.  First,
coproducts are unique up to canonical isomorphism.  Second, if we pick a
coproduct $x \sqcup y$ for every pair of objects $x,y \in \FinSet$, making
disjoint union into an operation in the traditional sense, we obtain
natural isomorphisms
\[   a_{x,y,z}\maps (x \sqcup y) \sqcup z \to x \sqcup (y \sqcup z),  \]
\[                l_x \maps \emptyset \sqcup x \to x, \qquad 
            r_x \maps x \sqcup \emptyset \to x .\]
Third, these natural isomorphisms satisfy coherence laws making 
$\FinSet$ into a monoidal category.   In short, the standard definition 
of monoidal category, which lists coherence laws, is best regarded as a
spinoff of the fact that monoidal structures defined by universal
properties automatically satisfy these laws.   

Studying this example also suggests another idea which is built into our
definition of $n$-categories.  To prove the existence of the associator
$a_{x,y,z}$ one uses the universal property for the coproduct of {\it
three} objects, and to prove the existence of $l_x$ and $r_x$ one
uses the universal property for the coproduct of {\it one} object.
This suggests that in an $n$-category, arbitrary $\ell$-ary composites
should be treated on an equal footing with binary composites.  The
formalism of operads is admirably suited for this task.

In what follows we first review the theory of typed operads.   Then we
use this to define the opetopes, and more generally, `$O$-opetopes' for
any typed operad $O$.   After a brief discussion  of some notions
concerning $O$-opetopic sets, we give the the definition of
$n$-categories, and more generally, $n$-coherent $O$-algebras.  We skim
over some technical details which can be found in our paper \cite{BD2}.

\subsection{Typed Operads}

To describe `many-sorted' or `typed' algebraic structures using operads,
we need to generalize the concept of operad a bit.   For any set $S$ of
`types', there is a notion of `$S$-operad'.   The basic idea is that for
any $x_1,\dots,x_\ell,x' \in S$, an $S$-operad $O$ has a set
$O(x_1,\dots,x_\ell;x')$ of  $k$-ary operations with inputs of type
$x_1,\dots,x_\ell$  and output of type $x'$.  As in an ordinary untyped
operad, we can visualize such an operation as a tree, but now we label
the edges of the tree by types.  For example, an operation $f \in
O(x_1,\dots,x_4;x')$ is drawn as follows:

\begin{center} \setlength{\unitlength}{0.000500in}
\begingroup\makeatletter\ifx\SetFigFont\undefined 
characters in \fmtname \def\x#1#2#3#4#5#6#7\relax{\def\x{#1#2#3#4#5#6}}%
\expandafter\x\fmtname xxxxxx\relax \def\y{splain}
or SliTeX? \gdef\SetFigFont#1#2#3{%
  \ifnum #1<17\tiny\else \ifnum #1<20\small\else
  \ifnum #1<24\normalsize\else \ifnum #1<29\large\else
  \ifnum #1<34\Large\else \ifnum #1<41\LARGE\else
     \huge\fi\fi\fi\fi\fi\fi
  \csname #3\endcsname}%
\else
\gdef\SetFigFont#1#2#3{\begingroup
  \count@#1\relax \ifnum 25<\count@\count@25\fi
  \def\x{\endgroup\@setsize\SetFigFont{#2pt}}%
  \expandafter\x
    \csname \romannumeral\the\count@ pt\expandafter\endcsname
    \csname @\romannumeral\the\count@ pt\endcsname
  \csname #3\endcsname}%
\fi\endgroup
\begin{picture}(1897,2565)(3226,-3304)
\thicklines
\put(4201,-1861){\circle*{150}}
\put(4201,-3061){\line( 0, 1){1200}}
\put(4201,-1861){\line(-1, 3){300}}
\put(4201,-1861){\line( 1, 3){300}}
\put(4201,-1861){\line( 1, 1){900}}
\put(4201,-1861){\line(-1, 1){900}}
\put(3226,-811){\makebox(0,0)[lb]{\raisebox{0pt}[0pt][0pt]{$x_1$}}}
\put(4126,-3286){\makebox(0,0)[lb]{\raisebox{0pt}[0pt][0pt]{$x'$}}}
\put(3901,-811){\makebox(0,0)[lb]{\raisebox{0pt}[0pt][0pt]{$x_2$}}}
\put(4501,-811){\makebox(0,0)[lb]{\raisebox{0pt}[0pt][0pt]{$x_3$}}}
\put(5101,-811){\makebox(0,0)[lb]{\raisebox{0pt}[0pt][0pt]{$x_4$}}}
\put(4426,-1966){\makebox(0,0)[lb]{\raisebox{0pt}[0pt][0pt]{$f$}}}
\end{picture}
\end{center}
\medskip

More precisely:
\begin{defn}\et For any set $S$, an `$S$-operad' $O$ consists of
{\rm \begin{enumerate}  
\item {\it for any $x_1,\dots,x_\ell,x' \in S$, a set
$O(x_1,\dots,x_\ell;x')$}
\item  {\it for any $f \in O(x_1,\dots,x_\ell;x')$ and any
$g_1 \in O(x_{11},\dots,x_{1i_1};x_1),$
$\dots, g_\ell \in O(x_{\ell 1},\dots,x_{\ell i_\ell};x_\ell)$, an element
\[ f \cdot (g_1,\dots, g_\ell) \in O(x_{11},\dots,x_{1i_1}, \dots \dots,
x_{\ell 1}, \dots , x_{\ell i_\ell};x') \] }
\item {\it for any $x \in S$, an element $1_x \in O(x;x)$}
\item {\it for any permutation $\sigma \in S_\ell$, a map
\ban     \sigma \maps O(x_1,\dots,x_\ell;x') &\to&
          O(x_{\sigma(1)}, \dots, x_{\sigma(\ell)};x')  \\
f &\mapsto& f\sigma \ean  }
\end{enumerate}}
\noindent such that:
\begin{alphalist}
\item  whenever both sides make sense, 
\[      f \cdot (g_1 \cdot (h_{11}, \dots, h_{1i_1}), \dots, 
g_\ell \cdot (h_{\ell 1}, \dots, h_{\ell i_\ell})) = \]
\[ (f \cdot (g_1, \dots g_\ell)) \; \cdot \; (h_{11}, \dots, h_{1i_1},
\dots\dots , h_{\ell 1}, \dots, h_{\ell i_\ell}) \]
\item  for any $f \in O(x_1,\dots,x_\ell;x')$,
\[ f =  1_{x'} \cdot f = f \cdot (1_{x_1},\dots,1_{x_\ell}) \]
\item  for any $f \in O(x_1,\dots,x_\ell;x')$ and 
$\sigma,\sigma' \in S_\ell$, 
\[      f(\sigma \sigma') = (f \sigma)\sigma' \]
\item  for any $f \in O(x_1,\dots,x_\ell;x')$, $\sigma \in S_\ell$, and 
$g_1 \in O(x_{11},\dots,x_{1i_1};x_1), $ \hfill \break
$\dots, g_\ell \in O(x_{\ell 1},\dots,x_{\ell i_\ell};x_\ell)$, 
\[     (f\sigma) \cdot (g_{\sigma(1)}, \dots, g_{\sigma(\ell)}) = (f \cdot
(g_1, \dots, g_\ell))\, \rho(\sigma), \]
where $\rho \maps S_\ell \to S_{i_1 + \cdots + i_\ell}$ is the obvious
homomorphism.
\item  for any $f \in O(x_1,\dots,x_\ell;x')$, 
$g_1 \in O(x_{11},\dots,x_{1i_1};x_1),\dots,$ \hfill \break
$g_\ell \in O(x_{k1},\dots,x_{\ell i_\ell};x_\ell)$, and 
$\sigma_1 \in S_{i_1}, \dots, \sigma_\ell \in S_{i_\ell}$,
\[     (f\cdot (g_1\sigma_1 , \dots, g_\ell \sigma_\ell))  = (f \cdot
(g_1, \dots, g_\ell))\, \rho'(\sigma_1,\dots,\sigma_\ell), \]
where $\rho' \maps S_{i_1} \times \cdots \times S_{i_\ell}
\to S_{i_1 + \cdots + i_\ell}$ is the obvious homomorphism.
\end{alphalist}
\end{defn}

There is an obvious notion of a morphism from an $S$-operad $O$ to an
$S$-operad $O'$: a function mapping each set $O(x_1,\dots,x_\ell;x')$ to
the corresponding set $O'(x_1,\dots,x_\ell;x')$, preserving composition,
identities, and the symmetric group actions.  An important example is an
`algebra' of an $S$-operad:

\begin{defn}  \label{c.op.algebras}\et  For any $S$-operad
$O$, an `$O$-algebra' $A$ consists of:
{\rm \begin{enumerate}
\item {\it for any $x \in S$, a set $A(x)$.}
\item {\it for any $f \in O(x_1,\dots,x_\ell;x')$, a function
\[   \alpha(f) \maps A(x_1) \times \cdots \times A(x_\ell) \to A(x') \]  }
\end{enumerate}}
\noindent such that:
\begin{alphalist}
\item whenever both sides make sense,
\[      \alpha(f \cdot (g_1, \dots, g_\ell)) = \alpha(f) (\alpha(g_1)
\times \cdots \times \alpha(g_\ell)) \]
\item  for any $x \in C$, $\alpha(1_x)$ acts as the
identity on $A(x)$
\item for any $f \in O(x_1,\dots,x_\ell,x')$ and $\sigma \in S_\ell$,
\[         \alpha(f\sigma) = \alpha(f)\sigma, \]
where $\sigma \in S_\ell$ acts on the function $\alpha(f)$ 
on the right by permuting its arguments. 
\end{alphalist}
\end{defn}

In what follows, by `operad' we will mean an $S$-operad for some set $S$
of types.  We can think of such an operad as a simple sort of theory,
and its algebras as models of this theory.   Thus we can study operads
either `syntactically' or `semantically'.  To describe an operad
syntactically, we list: 
\begin{enumerate} 
\item the set $S$ of {\it types}, 
\item the sets $O(x_1,\dots,x_\ell;x')$ of {\it operations}, 
\item the set of all {\it reduction laws} saying that some composite of
operations (possibly with arguments permuted) equals some other
operation.   
\end{enumerate} 
This is like a presentation in terms of generators and relations, with
the reduction laws playing the role of relations.  On the other hand, to
describe an operad semantically, we describe its algebras.  

For example, the simplest operad is the `initial untyped operad' $I$.
Syntactically, this is the $S$-operad with:
\begin{enumerate}
\item only one type: $S = \{x\}$,
\item only one operation, the identity operation $1 \in O(x;x)$,
\item all possible reduction laws.
\end{enumerate}
Semantically, $I$ is the operad whose algebras are just sets.

Another important operad is the `terminal untyped operad' $T$.  
This is the $S$-operad with 
\begin{enumerate}
\item only one type: $S = \{x\}$,
\item exactly one operation of each arity,
\item all possible reduction laws.
\end{enumerate}
The algebras of $T$ are commutative monoids, with the $\ell$-ary
operation being $\ell$-fold multiplication, or the unit element when
$\ell = 0$, since nullary operations correspond to `constants'.

\subsection{Opetopes}

The following fact is the key to defining the opetopes.  Let $O$ be an
$S$-operad, and let $\elt(O)$ be the set of all operations of $O$.

\begin{thm}\et There is an $\elt(O)$-operad $O^+$ whose algebras are
$S$-operads over $O$, i.e., $S$-operads equipped with a homomorphism to
$O$.  \end{thm}

\noindent We call $O^+$ the `slice operad' of $O$.  One can describe 
$O^+$ syntactically as follows:
\begin{enumerate}
\item The types of $O^+$ are the operations of $O$.
\item The operations of $O^+$ are the reduction laws of $O$.
\item The reduction laws of $O^+$ are the ways of combining reduction
laws of $O$ to give other reduction laws.
\end{enumerate}
The `level-shifting' going on here as we pass from $O$ to $O^+$ is a way
of systematizing the process of categorification.  

A nice example of the slice operad construction is the operad $I^+$.
The algebras of this operad turn out to be monoids.  Another nice
example is the operad $T^+$, whose algebras are untyped operads!
It is also very interesting to iterate the slice operad construction.  
For example, the algebras of $I^{++}$ work out to be exactly untyped 
planar operads.  More generally, let $O^{n+}$ be the operad formed 
by applying the slice operad construction $n$ times to the operad $O$, 
or just $O$ itself if $n = 0$.   

\begin{defn}\et An $n$-dimensional `$O$-opetope' is a type of $O^{n+}$,
or equivalently, if $n \ge 1$, an operation of $O^{(n-1)+}$.  \end{defn}

In particular, we define an $n$-dimensional `opetope' to be an
$n$-dimensional $O$-opetope for $O = I$, the initial untyped operad.
The 0-dimensional opetopes are thus the types of $I$, but there is only
one type, so there is only one 0-dimensional opetope, which we visualize
as a point. The 1-dimensional opetopes are the types of $I^+$, or in
other words, the operations of $I$.  $I$ has only one operation, the
identity, so there is only one 1-dimensional opetope, which we visualize
as an interval.  The 2-dimensional opetopes are the types of $I^{++}$,
or in other words, the operations of $I^+$, which are the reduction laws
of $I$.  These reduction laws all state that the identity operation
composed with itself $\ell$ times equals itself.  This leads to
2-dimensional opetopes with $\ell$ infaces and one outface.  Actually
there are $\ell !$ different 2-dimensional opetopes with $\ell$ infaces,
since the permutation group $S_\ell$ acts freely on the set of
$\ell$-ary operations of $I^+$.  We could keep track of these by
labelling the infaces with some permutation of $\ell$ distinct symbols.
A more systematic approach is to use `metatree notation'.  In this notation, 
any $n$-dimensional $O$-opetope is represented as a list of $n$ labelled
trees.  

To see how this works, first consider the 2-dimensional $O$-opetopes, which
are the operations of $O^+$.   An operation of $O^+$ can be specified as
in the figure below.  

\medskip
\vbox{\begin{center}

\setlength{\unitlength}{0.000500in}%
\begingroup\makeatletter\ifx\SetFigFont\undefined
\def\x#1#2#3#4#5#6#7\relax{\def\x{#1#2#3#4#5#6}}%
\expandafter\x\fmtname xxxxxx\relax \def\y{splain}%
\ifx\x\y   
\gdef\SetFigFont#1#2#3{%
  \ifnum #1<17\tiny\else \ifnum #1<20\small\else
  \ifnum #1<24\normalsize\else \ifnum #1<29\large\else
  \ifnum #1<34\Large\else \ifnum #1<41\LARGE\else
     \huge\fi\fi\fi\fi\fi\fi
  \csname #3\endcsname}%
\else
\gdef\SetFigFont#1#2#3{\begingroup
  \count@#1\relax \ifnum 25<\count@\count@25\fi
  \def\x{\endgroup\@setsize\SetFigFont{#2pt}}%
  \expandafter\x
    \csname \romannumeral\the\count@ pt\expandafter\endcsname
    \csname @\romannumeral\the\count@ pt\endcsname
  \csname #3\endcsname}%
\fi
\fi\endgroup
\begin{picture}(8744,5144)(579,-5183)
\thicklines
\put(2401,-3661){\circle*{150}}
\put(3301,-2761){\circle*{150}}
\put(1501,-2761){\circle*{150}}
\put(901,-1561){\circle*{150}}
\put(7801,-3631){\circle*{150}}
\put(2101,-1561){\circle*{150}}
\put(2401,-2761){\circle*{150}}
\put(2401,-3661){\line( 0,-1){1500}}
\put(2401,-3661){\line(-1, 1){900}}
\put(2401,-3661){\line( 1, 1){900}}
\put(1501,-2761){\line(-1, 2){600}}
\put(1501,-2761){\line( 1, 2){600}}
\put(7801,-5161){\line( 0, 1){1500}}
\put(901,-1561){\line( 6, 5){1800}}
\put(3301,-2686){\line( 1, 4){652.941}}
\put(2101,-1561){\line(-1, 1){450}}
\put(1426,-886){\line(-1, 1){825}}
\put(2401,-3661){\line( 0, 1){900}}
\put(7801,-3661){\line(-1, 5){300}}
\put(7801,-3661){\line( 1, 5){300}}
\put(7801,-3661){\line(-3, 5){900}}
\put(7801,-3661){\line( 3, 5){900}}
\put(7801,-3661){\line(-1, 1){1500}}
\put(7801,-3661){\line( 1, 1){1500}}
\put(601,-1786){\makebox(0,0)[lb]{\raisebox{0pt}[0pt][0pt]{$f_1$}}}
\put(1201,-3061){\makebox(0,0)[lb]{\raisebox{0pt}[0pt][0pt]{$f_4$}}}
\put(2026,-3811){\makebox(0,0)[lb]{\raisebox{0pt}[0pt][0pt]{$f_3$}}}
\put(3451,-2986){\makebox(0,0)[lb]{\raisebox{0pt}[0pt][0pt]{$B$}}}
\put(1726,-1711){\makebox(0,0)[lb]{\raisebox{0pt}[0pt][0pt]{$f_5$}}}
\put(2256,-1636){\makebox(0,0)[lb]{\raisebox{0pt}[0pt][0pt]{$E$}}}
\put(1696,-2836){\makebox(0,0)[lb]{\raisebox{0pt}[0pt][0pt]{$D$}}}
\put(1126,-1636){\makebox(0,0)[lb]{\raisebox{0pt}[0pt][0pt]{$A$}}}
\put(2626,-3736){\makebox(0,0)[lb]{\raisebox{0pt}[0pt][0pt]{$C$}}}
\put(2101,-2911){\makebox(0,0)[lb]{\raisebox{0pt}[0pt][0pt]{$f_6$}}}
\put(2451,-2836){\makebox(0,0)[lb]{\raisebox{0pt}[0pt][0pt]{$F$}}}
\put(2926,-2761){\makebox(0,0)[lb]{\raisebox{0pt}[0pt][0pt]{$f_2$}}}
\put(6226,-2086){\makebox(0,0)[lb]{\raisebox{0pt}[0pt][0pt]{$A$}}}
\put(6826,-2086){\makebox(0,0)[lb]{\raisebox{0pt}[0pt][0pt]{$B$}}}
\put(7426,-2086){\makebox(0,0)[lb]{\raisebox{0pt}[0pt][0pt]{$C$}}}
\put(8026,-2086){\makebox(0,0)[lb]{\raisebox{0pt}[0pt][0pt]{$D$}}}
\put(8626,-2086){\makebox(0,0)[lb]{\raisebox{0pt}[0pt][0pt]{$E$}}}
\put(9226,-2086){\makebox(0,0)[lb]{\raisebox{0pt}[0pt][0pt]{$F$}}}
\end{picture}
\end{center}

\medskip
}

\noindent The first tree is an arbitrary `$O$-tree'.  This is a
combed tree (i.e., planar except for a permutation of inputs at the top)
with nodes labelled by operations of $O$.  We require that a node
labelled by a $\ell$-ary operation have $\ell$ edges coming into it from
above.  Moreover, we require that it be possible to label every edge
with an element of $S$ in such a way that for any node labelled by an
operation in $O(x_1,\dots,x_\ell;x')$, the edges coming into that node
from above are labelled by the types $x_1,\dots,x_\ell$ in that order,
while the edge coming out of it from below is labelled by the type $x'$.
We also label each node of this first tree with a distinct symbol
$A,B,C,$ etc..  The second tree is planar and has only one node, with
$n$ edges coming into that node from above, labelled by the same symbols
$A,B,C,\dots$ in any order.  These specify the order of the input types
of the operation of $O^+$ we are describing.

More generally, for any $n > 1$ one can specify any $n$-dimensional
$O$-opetope by means of an `$n$-dimensional metatree'.  Here is an
example for $n = 3$:

\medskip
\vbox{
\begin{center}

\setlength{\unitlength}{0.00037500in}%
\begingroup\makeatletter\ifx\SetFigFont\undefined
\def\x#1#2#3#4#5#6#7\relax{\def\x{#1#2#3#4#5#6}}%
\expandafter\x\fmtname xxxxxx\relax \def\y{splain}%
\ifx\x\y   
\gdef\SetFigFont#1#2#3{%
  \ifnum #1<17\tiny\else \ifnum #1<20\small\else
  \ifnum #1<24\normalsize\else \ifnum #1<29\large\else
  \ifnum #1<34\Large\else \ifnum #1<41\LARGE\else
     \huge\fi\fi\fi\fi\fi\fi
  \csname #3\endcsname}%
\else
\gdef\SetFigFont#1#2#3{\begingroup
  \count@#1\relax \ifnum 25<\count@\count@25\fi
  \def\x{\endgroup\@setsize\SetFigFont{#2pt}}%
  \expandafter\x
    \csname \romannumeral\the\count@ pt\expandafter\endcsname
    \csname @\romannumeral\the\count@ pt\endcsname
  \csname #3\endcsname}%
\fi
\fi\endgroup
\begin{picture}(12944,5305)(579,-5183)
\thicklines
\put(2401,-3661){\circle*{150}}
\put(3301,-2761){\circle*{150}}
\put(1501,-2761){\circle*{150}}
\put(901,-1561){\circle*{150}}
\put(2101,-1561){\circle*{150}}
\put(2401,-2761){\circle*{150}}
\put(5101,-1861){\circle*{150}}
\put(5701,-2761){\circle*{150}}
\put(6601,-3661){\circle*{150}}
\put(10201,-3661){\circle*{150}}
\put(10801,-2461){\circle*{150}}
\put(12901,-3661){\circle*{150}}
\put(2401,-3661){\line( 0,-1){1500}}
\put(2401,-3661){\line(-1, 1){900}}
\put(2401,-3661){\line( 1, 1){900}}
\put(1501,-2761){\line(-1, 2){600}}
\put(1501,-2761){\line( 1, 2){600}}
\put(901,-1561){\line( 6, 5){1800}}
\put(3301,-2686){\line( 1, 4){652.941}}
\put(2101,-1561){\line(-1, 1){450}}
\put(1426,-886){\line(-1, 1){825}}
\put(2401,-3661){\line( 0, 1){900}}
\put(6601,-5161){\line( 0, 1){1500}}
\put(6601,-3661){\line( 1, 1){900}}
\put(6601,-3661){\line(-1, 1){900}}
\put(5701,-2761){\line(-2, 3){600}}
\put(5701,-2761){\line( 2, 3){600}}
\put(5701,-2761){\line( 0, 1){900}}
\put(5101,-1861){\line(-1, 4){300}}
\put(5101,-1861){\line( 0, 1){1200}}
\put(5101,-661){\line( 0, 1){600}}
\put(4801,-661){\line(-1, 4){150}}
\put(8101,-361){\line( 0, 1){300}}
\put(7501,-2761){\line( 1, 3){600}}
\put(8101,-961){\line( 0, 1){600}}
\put(5101,-1861){\line( 1, 3){300}}
\put(5401,-961){\line( 5, 3){1500}}
\put(5701,-1861){\line( 0, 1){975}}
\put(5701,-661){\line( 0, 1){600}}
\put(6301,-1861){\line( 0, 1){1350}}
\put(6301,-511){\line( 0,-1){ 75}}
\put(6301,-586){\line( 0, 1){ 75}}
\put(6301,-61){\line( 0,-1){225}}
\put(10201,-5161){\line( 0, 1){1500}}
\put(10201,-3661){\line( 1, 2){600}}
\put(10201,-3661){\line(-1, 2){1200}}
\put(10801,-2461){\line( 1, 2){600}}
\put(10801,-2536){\line(-1, 2){630}}
\put(12901,-5161){\line( 0, 1){1500}}
\put(12901,-3661){\line( 1, 2){600}}
\put(12901,-3661){\line(-1, 2){600}}
\put(9001,-1261){\line( 1, 1){1200}}

\put(10181,-1261){\line(-1, 1){525}}

\put(9001,-61){\line( 1,-1){525}}
\put(11401,-1261){\line( 0, 1){1200}}
\put(551,-1786){\makebox(0,0)[lb]{\raisebox{0pt}[0pt][0pt]{$f_1$}}}
\put(1151,-3061){\makebox(0,0)[lb]{\raisebox{0pt}[0pt][0pt]{$f_4$}}}
\put(1826,-3811){\makebox(0,0)[lb]{\raisebox{0pt}[0pt][0pt]{$f_3$}}}
\put(3451,-2986){\makebox(0,0)[lb]{\raisebox{0pt}[0pt][0pt]{$B$}}}
\put(1556,-1711){\makebox(0,0)[lb]{\raisebox{0pt}[0pt][0pt]{$f_5$}}}
\put(2326,-1636){\makebox(0,0)[lb]{\raisebox{0pt}[0pt][0pt]{$E$}}}
\put(1126,-1636){\makebox(0,0)[lb]{\raisebox{0pt}[0pt][0pt]{$A$}}}
\put(2626,-3736){\makebox(0,0)[lb]{\raisebox{0pt}[0pt][0pt]{$C$}}}
\put(1971,-2911){\makebox(0,0)[lb]{\raisebox{0pt}[0pt][0pt]{$f_6$}}}
\put(2816,-2761){\makebox(0,0)[lb]{\raisebox{0pt}[0pt][0pt]{$f_2$}}}
\put(2476,-2611){\makebox(0,0)[lb]{\raisebox{0pt}[0pt][0pt]{$F$}}}
\put(1651,-2761){\makebox(0,0)[lb]{\raisebox{0pt}[0pt][0pt]{$D$}}}
\put(4576, 14){\makebox(0,0)[lb]{\raisebox{0pt}[0pt][0pt]{$A$}}}
\put(5026, 14){\makebox(0,0)[lb]{\raisebox{0pt}[0pt][0pt]{$D$}}}
\put(5626, 14){\makebox(0,0)[lb]{\raisebox{0pt}[0pt][0pt]{$C$}}}
\put(6826, 14){\makebox(0,0)[lb]{\raisebox{0pt}[0pt][0pt]{$E$}}}
\put(6226, 14){\makebox(0,0)[lb]{\raisebox{0pt}[0pt][0pt]{$F$}}}
\put(8026, 14){\makebox(0,0)[lb]{\raisebox{0pt}[0pt][0pt]{$B$}}}
\put(4706,-2011){\makebox(0,0)[lb]{\raisebox{0pt}[0pt][0pt]{$G$}}}
\put(5306,-2986){\makebox(0,0)[lb]{\raisebox{0pt}[0pt][0pt]{$H$}}}
\put(6171,-3811){\makebox(0,0)[lb]{\raisebox{0pt}[0pt][0pt]{$I$}}}
\put(11026,-2611){\makebox(0,0)[lb]{\raisebox{0pt}[0pt][0pt]{$J$}}}
\put(10426,-3811){\makebox(0,0)[lb]{\raisebox{0pt}[0pt][0pt]{$K$}}}
\put(13426,-2386){\makebox(0,0)[lb]{\raisebox{0pt}[0pt][0pt]{$J$}}}
\put(12226,-2386){\makebox(0,0)[lb]{\raisebox{0pt}[0pt][0pt]{$K$}}}
\put(8926, 14){\makebox(0,0)[lb]{\raisebox{0pt}[0pt][0pt]{$G$}}}
\put(11326, 14){\makebox(0,0)[lb]{\raisebox{0pt}[0pt][0pt]{$H$}}}
\put(10126, 14){\makebox(0,0)[lb]{\raisebox{0pt}[0pt][0pt]{$I$}}}
\end{picture}
\end{center}

\medskip
}

An $n$-dimensional metatree is a list of $n$ labelled trees, the last of
which is a planar tree with only one node, while the rest are combed
trees. The first tree is an arbitrary $O$-tree.  For $1 \le i < n$,
every node of the $i$th tree is labelled with a distinct symbol, and the
same symbols also label all the edges at the very top of the $(i + 1)$st
tree, each symbol labelling exactly one edge.  In addition, each edge of
the $(i + 1)$st tree must correspond to a subtree of the $i$th tree in
such a way that:  
\begin{enumerate} 
\item The edge at the very top of
the $(i + 1)$st tree labelled by a given symbol corresponds to the
subtree of the $i$th tree whose one and only node is labelled by the
same symbol.   
\item  The edge of the $(i + 1)$st tree coming out of a
given note from below corresponds to the subtree that is the union of
the subtrees corresponding to the edges coming into that node from
above. 
\item  The edge at the very bottom of the $(i+1)$st tree
corresponds to the whole $i$th tree. 
\end{enumerate} 
Special care must be taken when the node of the last tree has no edges
coming into it from above.  This can only occur when all the previous
trees are empty.  This sort of metatree describes a nullary operation of
$O^{(n-1)+}$ whose output type is an identity operation $1_x$ of
$O^{(n-2)+}$.  To specify which identity operation, we need to label the
edge coming out of the node of the last tree from below with the
operation $1_x$.  

\section{$n$-Coherent operad algebras}

An `$n$-coherent $O$-algebra' is an $n$ times categorified analog of an
algebra of the operad $O$.  In particular, when $O = I$, an $n$-coherent
$O$-algebra is just an $n$-category, which is the $n$ times categorified
analog of a set.   An $n$-coherent $O$-algebra is an `$O$-opetopic set'
with certain properties.  We omit the precise definition of 
$O$-opetopic sets here.  For our purposes, it should suffice to know
that an $O$-opetopic set is very much like a simplicial set, but with
$O$-opetopes replacing simplices, and no `degeneracy maps', only `face
maps'.  An $O$-opetopic set thus consists of (possibly empty) 
collections of `cells' shaped like all the different $O$-opetopes, such
that any inface or outface of a cell is again a cell.   

If $j \ge 1$, we may schematically represent a $j$-dimensional
cell $x$ in an $O$-opetopic set as follows:
\[ 
\begin{diagram}[(a_1,\dots,a_\ell)]
\node{(a_1,\dots,a_\ell)}\arrow{e,t}{x} \node{a'}
\end{diagram} 
\]
Here $a_1,\dots,a_\ell$ are the infaces of $x$ and $a'$ is the outface
of $x$; all these are cells of one lower dimension.   A configuration
just like this, satisfying all the incidence relations satisfied by the
boundary of a cell, but with $x$ itself missing:
\[
\begin{diagram}[(a_1,\dots,a_\ell)]
\node{(a_1,\dots,a_\ell)}\arrow{e,t}{?} \node{a'}
\end{diagram}
\] 
is called a `frame'.  A `niche' is like a frame with the outface missing:
\[ 
\begin{diagram}[(a_1,\dots,a_\ell)]
\node{(a_1,\dots,a_\ell)}\arrow{e,t}{?} \node{?}
\end{diagram} 
\]
Similarly, a `punctured niche' is like a frame with the outface
and one inface missing:
\[ 
\begin{diagram}[(a_{1},a_{i+1},\dots,a_{\ell})]
\node{(a_{1},\dots,a_{i-1},?,a_{i+1},\dots,a_{\ell})}\arrow{e,t}{?} \node{?}
\end{diagram} 
\]
If one of these configurations (frame, niche, or punctured niche) can be
extended to an actual cell, the cell is called an `occupant' of the
configuration.  Occupants of the same frame are called
`frame-competitors', while occupants of the same niche are called
`niche-competitors'.  

Next we need the concept of a `universal occupant' of a niche.  Since
the definition of this concept looks rather formidable at first, we
first give a heuristic explanation.  As already noted, the main use
of universality to define composites:

\begin{defn}\et  Given a universal occupant $u$ of a $j$-dimensional niche:
\[
\begin{diagram}[(a_1,\dots,a_k)]
\node{(a_1,\dots,a_k)}\arrow{e,t}{u} \node{b}
\end{diagram}
\]
we call $b$ a `composite' of $(a_{1},\dots,a_{k})$.
\end{defn}

To understand universality more deeply, one must understand the
role played by cells of different dimensions.  In our framework an
$n$-category usually has cells of arbitrarily high dimension, just like
a Kan complex.  For $j \le n$ the $j$-dimensional cells play the role of
$j$-morphisms, while for $j > n$ they play the role of `equations',
`equations between equations', and so on.  The definition of
universality depends on $n$ in a way that has the following effects. 
For $j \le n$ there may be many universal occupants of a given
$j$-dimensional niche, which is why we speak of `a' composite rather
than `the' composite.  There is at most one occupant of any given
$(n+1)$-dimensional niche, which is automatically universal.  Thus
composites of $n$-cells are unique, and we may think of the universal
occupant of an $(n+1)$-dimensional niche as an equation saying that the
composite of the infaces equals the outface.  For $j > n + 1$ there is
exactly one occupant of each $j$-dimensional frame, indicating that the
composite of the equations corresponding to the infaces equals the
equation corresponding to the outface.

The definition of universality essentially says that a $j$-dimensional
niche-occupant is universal if all of its niche-competitors factor
through it uniquely, {\it up to equivalence}.  For $j \ge n+1$ this
amounts to saying that each niche has a unique occupant, while for $j =
n$ it means that each niche has an occupant through which all of its
niche-competitors factor uniquely.  Technically, the definition of
universality says that composition with a universal niche-occupant set
up a `balanced punctured niche' of one higher dimension.  One should
think of a balanced punctured niche as defining an equivalence between
occupants of its outface and occupants of its missing outface.

Now let us give the actual definition:

\begin{defn}\et A $j$-dimensional niche-occupant:
\[
\begin{diagram}[(c_1,\dots,c_k)]
\node{(c_1,\dots,c_k)}  \arrow{e,t}{u} \node{d}
\end{diagram}
\]
is said to be `universal' if and only if $j>n$ and $u$ is the only
occupant of its niche, 
or $j\le n$ and for any frame-competitor $d'$ of $d$, the
$(j+1)$-dimensional punctured niche
\[
\begin{diagram}[((c_{1},\dots,c_{k}) \to d,\; d \to d')]
\node{((c_{1},\dots,c_{k}) \mapright{u} d,\; d \mapright{?} d')} 
\arrow{s,r}{?} \\
\node{(c_{1},\dots,c_{k}) \mapright{?} d'}  
\end{diagram}
\]
and its mirror-image version
\[
\begin{diagram}[((c_{1},\dots,c_{k}) \to d,\; d \to d')]
\node{(d \mapright{?} d',\;(c_{1},\dots,c_{k}) \mapright{u} d)} 
\arrow{s,r}{?} \\
\node{(c_{1},\dots,c_{k}) \mapright{?} d'}  
\end{diagram}
\]
are balanced.  
\end{defn}

Of course, now we need the definition of `balanced'.  The reader will
note that while the definitions of `universal' and `balanced' call upon
each other, there is no bad circularity.

\begin{defn} \et An $m$-dimensional punctured niche:
\[ 
\begin{diagram}[(a_{1},a_{j+1},\dots,a_{k})]
\node{(a_{1},\dots,a_{i-1},?,a_{i+1},\dots,a_{k})}\arrow{e,t}{?} \node{?}
\end{diagram}
\]
is said to be `balanced' if and only if $m>n+1$ or:
\begin{enumerate}
\item any extension
\[ 
\begin{diagram}[(a_{1},a_{i+1},\dots,a_{m})]
\node{(a_{1},\dots,a_{i-1},?,a_{i+1},\dots,a_{k})}\arrow{e,t}{?} \node{b}
\end{diagram} 
\]
extends further to:
\[ 
\begin{diagram}[(a_{1},a_{i+1},\dots,a_{k})]
\node{(a_{1},\dots,a_{i-1},a_i,a_{i+1},\dots,a_{k})}\arrow{e,t}{u} \node{b}
\end{diagram} 
\]
with $u$ universal in its niche, and 
\item  for any occupant
\[ 
\begin{diagram}[(a_{1},a_{i+1},\dots,a_{k})]
\node{(a_{1},\dots,a_{i-1},a_i,a_{i+1},\dots,a_{k})}\arrow{e,t}{u} \node{b}
\end{diagram} 
\]
universal in its niche, and frame-competitor $a'_{i}$ of $a_{i}$, the
$(m+1)$-dimensional punctured niche
\[
\begin{diagram}
[((c_{1},\dots,c_{m}) \to d,\; d \to d')]
\node{(a'_{i} \mapright{?} a_{i}, \;
(a_{1},\dots,a_{i-1},a_{i},a_{i+1},\dots,a_{k}) \mapright{u} b)}  
\arrow{s,r}{?}  \\
\node{(a_{1},\dots,a_{i-1},a'_{i},a_{i+1},\dots, a_{k}) \mapright{?} b}
\end{diagram}
\]
and its mirror-image version
\[
\begin{diagram}
[((c_{1},\dots,c_{m}) \to d,\; d \to d')]
\node{((a_{1},\dots,a_{i-1},a_{i},a_{i+1},\dots,a_{k}) \mapright{u} b, \;
a'_{i} \mapright{?} a_{i})} 
\arrow{s,r}{?}  \\
\node{(a_{1},\dots,a_{i-1},a'_{i},a_{i+1},\dots, a_{k}) \mapright{?} b}
\end{diagram}
\]
are balanced.  
\end{enumerate}
\end{defn}

\noindent Note that the first numbered condition in the definition of
`balanced' definition generalizes the concept of an essentially
surjective functor, while the second generalizes the concept of a fully
faithful functor.

Finally, we define $n$-coherent $O$-algebras and various special cases:

\begin{defn} \et An `$n$-coherent $O$-algebra' is an $O$-opetopic
set such that 1) every niche has a universal occupant, and 2) composites
of universal cells are universal.  \end{defn}

\begin{defn} \et An `$n$-category' is an $n$-coherent $I$-algebra.
\end{defn}

\begin{defn} \et A `monoidal $n$-category' is an $n$-coherent $I^+$-algebra.
\end{defn}

\begin{defn} \et A `stable $n$-category' is an $n$-coherent $T$-algebra.
\end{defn}

One can show that any $S$-operad morphism $f \maps O \to O'$ allows one
to turn an $n$-coherent $O'$-algebra $A$ into an $n$-coherent
$O$-algebra $f^\ast A$.  Thus any stable $n$-category has an underlying
monoidal $n$-category, and any monoidal $n$-category has an underlying
$n$-category.  

Given an $n$-category with 0-cells $x$ and $y$, there is
an $(n-1)$-category $\hom(x,y)$.  One can also 
construct a stable $(n+1)$-category of all (small) $n$-categories,
though the details of this construction have not yet been published.
Using these facts, one can give rigorous formulations of many 
processes going between neighboring entries in Table 2:
decategorification, discrete categorification, delooping, looping,
forgetting monoidal structure, stabilization, and the generalized center
construction.  However, it remains to make precise and prove the
Stabilization Hypothesis in this framework.  Basically, one wishes to
show that for $k \ge n+2$, the $(n+1)$-category of all stable
$n$-categories is equivalent to the full sub-$(n+1)$-category of all
$n$-categories having only one 0-cell and only one $j$-cell in each
frame for $0 < j < k$.  

\section{Conclusions} \label{conclusions}

In this paper we have discussed iterated categorifications and
stabilizations of some of the very simplest algebraic structures: the
natural numbers and the integers.  However, one can also categorify many
other concepts: vector spaces \cite{KV} and Hilbert spaces \cite{B},
group algebras \cite{CY}, algebras of formal power series
\cite{BD2,Joyal} and other Hopf algebras \cite{CF,CY2}, sheaves
\cite{Breen,Brylinski}, and so on.  Interesting results about these
familiar structures typically have interesting categorified analogs.  It
is clear, therefore, that the set-based mathematics we know and love is
just the tip of an immense iceberg of $n$-categorical, and ultimately
$\omega$-categorical, mathematics.

The prospect of exploring this huge body of new mathematics is both
exhilarating and daunting.  The basic philosophy is simple: {\it never
mistake equivalence for equality}.  The technical details, however, are
not so simple --- at least not yet.  To proceed efficiently it is
crucial that we gain a clearer understanding of the foundations before
rushing ahead with complicated constructions.  

Many basic questions remain open.  For example, how significant is the
fact that operads play a role both in the theory of $E_k$ spaces and the
definition of $n$-categories described above?  Operads are very
versatile, so this might at first seem to be a coincidence.  However,
there are deep relationships between operads, categorification, and the
theory of algebraic structures satisfying laws `up to coherent homotopy'
\cite{Stasheff2}.  In particular, Trimble \cite{Trimble} has pointed out an
interesting connection.  For many purposes it is best to think of $E_k$
spaces as algebras, not of the little $k$-cubes operad, but of a closely
related operad $F(k)$ discovered by Getzler and Jones \cite{GJ}.  
The space $F(k)_\ell$ is the Fulton-MacPherson compactification of 
the configuration space of $\ell$ points in $\R^k$ modulo translations 
and dilations \cite{FM}.  In particular, $F(1)_\ell$ is just 
$K_\ell \times S_\ell$, where $K_\ell$ is the $(\ell - 2)$-dimensional 
associahedron.  

Since homotopy $n$-types of $E_1$ spaces 
correspond to monoidal $n$-groupoids, while $n$-coherent $I^+$-algebras
are monoidal $n$-categories, one might expect a relationship between
$F(1)$ and the $I^+$-opetopes.  The associahedron $K_\ell$ has 
a cell decomposition having cells in one-to-one correspondence with
planar trees with $\ell$ leaves for which all nodes have at least one
edge coming in from above.  It follows that the cells in the corresponding
decomposition of $F(1)$ correspond to a certain class of 2-dimensional
$I^+$-opetopes, or equivalently, 3-dimensional opetopes.  Is there a deeper 
relation between opetopes and the associahedron?  This might shed new 
light on the origin of coherence laws.

\subsection*{Acknowledgments} We thank Michael Batanin, Dan Christensen,
Mark Hovey, Andr\'e Joyal, Carlos Simpson, Charles Rezk, Ross Street,
and Todd Trimble for helpful discussions and correspondence, and thank
Ezra Getzler and Mikhail Kapranov both for useful correspondence and for
organizing the conference at which this paper was presented.


{\small 

\begin{thebibliography}{10}

\bibitem{Adams} J.\ Adams, {\sl Infinite Loop Spaces,} Princeton
U.\ Press, Princeton, 1978.

\bibitem{B} J.\ Baez, Higher-dimensional algebra II: 2-Hilbert spaces,
{\sl Adv.\ Math.} {\bf 127} (1997), 125-189.

\bibitem{B2} J.\ Baez, An introduction to $n$-categories, {\sl 
7th Conference on Category Theory and Computer Science,} eds.\ 
E.\ Moggi and G.\ Rosolini, Lecture Notes in Computer
Science {\bf 1290}, Springer Verlag, Berlin, 1997.

\bibitem{BD} J.\ Baez and J.\ Dolan, Higher-dimensional algebra and
topological quantum field theory, {\sl Jour.\ Math.\ Phys. }{\bf 36}
(1995), 6073-6105.

\bibitem{BD2} J.\ Baez and J.\ Dolan, Higher-dimensional algebra III:
$n$-Categories and the algebra of opetopes, to appear in {\sl Adv.\ Math.},
preprint available as q-alg/9702014.

\bibitem{BL} J.\ Baez and L.\ Langford, 2-Tangles, to appear in  {\sl
Lett.\ Math.\ Phys.}, preprint available online as q-alg/9703033.
 
\bibitem{BL2} J.\ Baez and L.\ Langford, Higher-dimensional algebra
IV: 2-Tangles, to appear in {\sl Adv.\ Math.}.

\bibitem{BN} J.\ Baez and M.\ Neuchl, Higher-dimensional algebra I:
Braided monoidal 2-categories, {\sl Adv. Math.\ }{\bf 121} (1996),
196-244.

\bibitem{BP} M.\ G.\ Barratt and S.\ Priddy, On the homology
of non-connected monoids and their associated groups, {\sl Comm.\ 
Math.\ Helv.} {\bf 47} (1972), 1-14.

\bibitem{Batanin} M.\ Batanin, Monoidal globular categories as a
natural environment for the theory of weak $n$-categories, to 
appear in {\sl Adv.\ Math.}, preprint available at
http://www-math.mpce.mq.edu.au/$\sim$mbatanin/papers.html

\bibitem{Baues} H.-J.\ Baues, Homotopy types, in {\sl Handbook of
Algebraic Topology}, ed.\ I.\ M.\ James, Elsevier, New York, 1995.

\bibitem{Benabou} J.\ B\'enabou, Introduction to bicategories, Lecture
Notes in Mathematics {\bf 47}, Springer Verlag, Berlin, 1967, pp.\ 1-77.

\bibitem{BV} J.\ M.\ Boardman and R.\ M.\ Vogt, Homotopy-everything
H-spaces, {\sl Bull.\ Amer.\ Math.\ Soc.} {\bf 74} (1968), 1117-1122.

\bibitem{BV2} J.\ M.\ Boardman and R.\ M.\ Vogt, {\sl Homotopy invariant
structures on topological spaces}, Lecture Notes in Mathematics
{\bf 347}, Springer Verlag, Berlin, 1973.

\bibitem{Breen} L.\ Breen, On the classification of 2-gerbes and 2-stacks,
{\sl Ast\'erisque} {\bf 225} (1994), 1-160.

\bibitem{Brown} R.\ Brown, From groups to groupoids: a brief survey, 
{\sl Bull.\ London Math.\ Soc.\ }{\bf 19} (1987), 113-134.

\bibitem{Brylinski} J.-L.\ Brylinski, {\sl Loop Spaces, Characteristic
Classes, and Geometric Quantization,} Birkhauser, Boston, 1993.

\bibitem{BS} R.\ Brown and C.\ Spencer, $\cal G$-groupoids, crossed
modules, and the classifying space of a topological group, {\sl Proc.\ Kon.\
Akad.\ v.\ Wet.\ }{\bf 79} (1976), 296-302.

\bibitem{Crane} L.\ Crane, Clock and category: is quantum gravity
algebraic?, {\sl Jour.\ Math.\ Phys.\ }{\bf 36} (1995), 6180-6193.

\bibitem{CF} L.\ Crane and I.\ Frenkel, Four dimensional topological
quantum field theory, Hopf categories, and the canonical bases, {\sl
Jour.\ Math.\ Phys.\ }{\bf 35} (1994), 5136-5154.  

\bibitem{CY} L.\ Crane and D.\ Yetter, Examples of categorification, 
preprint available as q-alg/9607028.  

\bibitem{CY2} L.\ Crane and D.\ Yetter, Deformations of (bi)tensor
categories, preprint available as q-alg/9612011.

\bibitem{Crans} S.\ E.\ Crans, On combinatorial models for
higher dimensional homotopies, Ph.D.\ thesis, University of Utrecht,
Utrecht, 1991.

\bibitem{Crans2} S.\ E.\ Crans, Generalized centers of braided and sylleptic
monoidal 2-categories, to appear in {\sl Adv.\ Math.}, preprint
available at http://www-math.mpce.mq.edu.au/$\sim$crans/papers/papers.html

\bibitem{CRS} J.\ S.\ Carter, J.\ H.\ Rieger and M.\ Saito,  
A combinatorial description of knotted surfaces and their isotopies, 
{\sl Adv.\ Math.} {\bf 127} (1997), 1-51.

\bibitem{DS} B.\ Day and R.\ Street, Monoidal bicategories and Hopf
algebroids, {\sl Adv.\ Math.} {\bf 129} (1997), 99-157.

\bibitem{EH} B.\ Eckmann and P.\ Hilton, Group-like structures in
categories, {\sl Math.\ Ann.\ }{\bf 145} (1962), 227-255.

\bibitem{EK} S.\ Eilenberg and G.\ M.\ Kelly, Closed categories, in
{\sl Proceedings of the Conference on Categorical Algebra,} eds.\ S.\
Eilenberg {\it et al}, Springer Verlag, New York, 1966.

\bibitem{EKMM} A.\ D.\ Elmendorf, I.\ Kriz, M.\ Mandell and J.\ P.\ 
May, Modern foundations for stable homotopy theory, in {\sl Handbook of
Algebraic Topology}, ed.\ I.\ M.\ James, Elsevier, New York, 1995.

\bibitem{FY} P.\ Freyd and D.\ Yetter, Braided compact monoidal
categories with applications to low dimensional topology, {\sl
Adv.\ Math.\ }{\bf 77} (1989), 156-182.

\bibitem{FM} W.\ Fulton and R.\ MacPherson, A compactification of
configuration spaces, {\sl Ann.\ Math.\ }{\bf 139} (1994), 183-225.

\bibitem{GZ} P.\ Gabriel and M.\ Zisman, {\sl Calculus of Fractions
and Homotopy Theory}, Ergebnisse der Math.\ und ihrer Grenzgebiete,
vol.\ 35, Springer Verlag, Berlin, 1967.

\bibitem{GJ} E.\ Getzler and J.\ Jones, Operads, homotopy algebra, and
iterated integrals for double loop spaces, preprint available as
hep-th/9403055.

\bibitem{GPS} R.\ Gordon, A.\ J.\ Power, and R.\ Street, Coherence for
tricategories, {\sl Memoirs Amer.\ Math.\ Soc.} {\bf 117} (1995)
Number 558.

\bibitem{Gro} A.\ Grothendieck, Pursuing stacks, unpublished
manuscript, 1983, distributed from UCNW, Bangor, United Kingdom.

\bibitem{HMP} C.\ Hermida, M.\ Makkai, and J.\ Power, On weak
higher dimensional categories, preprint available at
http://hypatia.dcs.qmw.ac.uk/authors/M/ \break \hfill 
MakkaiM/papers/multitopicsets/ 

\bibitem{Joyal} A.\ Joyal, Une th\'eorie combinatoire des s\'eries 
formelles, {\sl Adv.\ Math.\ }{\bf 42} (1981), 1-82.

\bibitem{JS} A.\ Joyal and R.\ Street, Tortile Yang-Baxter operators in
tensor categories, {\sl J.\ Pure Appl.\ Alg.\ }{\bf 71} (1991), 43-51.

\bibitem{JS2} A.\ Joyal and R.\ Street, Braided tensor categories,
{\sl Adv.\ Math.\ }{\bf 102} (1993), 20-78.

\bibitem{KV} M.\ Kapranov and V.\ Voevodsky, 2-Categories and
Zamolodchikov tetrahedra equations, in {\sl Proc.\ Symp.\ Pure
Math.\ }{\bf 56} Part 2 (1994), AMS, Providence, pp.\ 177-260.

Braided monoidal 2-categories, 2-vector-spaces, and Zamolodchikov
tetrahedra equations, Northwestern U.\ preprint.

\bibitem{Kelly} G.\ Kelly, Coherence theorems for lax algebras
and distributive laws, Lecture Notes in Mathematics
{\bf 420}, Springer Verlag, Berlin, 1974, pp.\ 281-375.

\bibitem{KS} G.\ Kelly and R.\ Street, Review of the elements of
2-categories, Lecture Notes in Mathematics {\bf 420}, Springer Verlag, Berlin,
1974, pp.\ 75-103.

\bibitem{KT} V.\ Kharlamov and V.\ Turaev, On the definition of the
2-category of 2-knots, {\sl Amer.\ Math.\ Soc.\ Transl.\ }{\bf 174}
(1996), 205-221.

\bibitem{Laplaza} M.\ Laplaza, Coherence for distributivity,
Lecture Notes in Mathematics {\bf 281}, Springer Verlag, Berlin, 1972,
pp.\ 29-72.

\bibitem{MacLane} S.\ Mac Lane, Natural associativity and
commutativity, {\sl Rice U.\ Studies} {\bf 49} (1963), 28-46.

\bibitem{MacLane2} S.\ Mac Lane, {\sl Categories for the Working
Mathematician}, Springer Verlag, New York, 1971.

\bibitem{MW} S.\ Mac Lane and J.\ Whitehead, On the 3-type of a complex,
{\sl Proc.\ Nat.\ Acad.\ Sci.\ }{\bf 36} (1950), 41-58.

\bibitem{Majid} S.\ Majid, Representations, duals and doubles of
monoidal categories,  {\sl Suppl.\ Rend.\ Circ.\ Mat.\ Palermo,} Series
II {\bf 26} (1991), 197-206.  

\bibitem{May} J.\ P.\ May, {\sl Simplicial Objects in Algebraic Topology},
Van Nostrand, Princeton, 1968.

\bibitem{May2} J.\ P.\ May, {\sl The Geometry of Iterated Loop Spaces},
Lecture Notes in Mathematics {\bf 271}, Springer Verlag, Berlin, 1972.

\bibitem{MQRT} J.\ P.\ May, F.\ Quinn, N.\ Ray and J.\ Tornehave, 
{\sl $E_\infty$ Ring Spaces and $E_\infty$ Ring Spectra},
Lecture Notes in Mathematics {\bf 577}, Springer Verlag, Berlin, 1977.

\bibitem{Power} A.\ J.\ Power, Why tricategories?, {\it Info.\ Comp.\
}{\bf 120} (1995), 251-262.

\bibitem{Quillen} D.\ Quillen, {\sl Homotopical Algebra}, Lecture
Notes in Mathematics {\bf 43}, Springer Verlag, Berlin, 1967.

\bibitem{Rota} G.-C.\ Rota, {\sl Gian-Carlo Rota on Combinatorics: 
Introductory Papers and Commentaries,} ed.\ J.\ P.\ S. Kung, Birkhauser,
Boston, 1995.

\bibitem{Segal} G.\ Segal, Configuration spaces and iterated loop-spaces, 
{\sl Invent.\ Math.\ }{\bf 21} (1973), 213-222.

\bibitem{Shum} M.-C.\ Shum, Tortile tensor categories, 
{\sl Jour.\ Pure Appl.\ Alg.\ }{\bf 93} (1994), 57-110.

\bibitem{Simpson} C.\ Simpson, A closed model structure for
$n$-categories, internal $Hom$, $n$-stacks and generalized Seifert-Van
Kampen, preprint available as alg-geom/9704006.

\bibitem{Simpson2} C.\ Simpson, Limits in $n$-categories, preprint
available as alg-geom/9708010.

\bibitem{Stasheff} J.\ D.\ Stasheff, Homotopy associativity of
H-spaces I \& II, {\sl Trans.\ Amer.\ Math.\ Soc.\ }{\bf 108} (1963),
275-292, 293-312. 

\bibitem{Stasheff2} J.-L.\ Loday, J.\ D.\ Stasheff, and A.\ A.\
Voronov, eds., {\sl Operads: Proceedings of Renaissance Conferences,}
A.\ M.\ S., Providence, Rhode Island, 1997.

\bibitem{Street} R.\ Street, The algebra of oriented simplexes, 
{\sl Jour.\ Pure Appl.\ Alg.\ }{\bf 49} (1987), 283-335.

\bibitem{Tamsamani} Z.\ Tamsamani, Sur des notions de
$\infty$-categorie et $\infty$-groupoide non-strictes via des
ensembles multi-simpliciaux, preprint available as
alg-geom/9512006.

\bibitem{Tamsamani2} Z.\ Tamsamani, Equivalence de la th\'eorie
homotopique des $n$-groupoides et celle des espaces 
topologiques $n$-tronqu\'es, preprint available as 
alg-geom/9607010.

\bibitem{Trimble} T.\ Trimble, personal communication. 

\bibitem{T} V.\ Turaev, Operator invariants of tangles, and
R-matrices, {\sl Math.\ USSR Izvestia} {\bf 35} (1990), 411-444.

\bibitem{Y} D.\ Yetter, Markov algebras, in {\sl Braids}, eds.\
J.\ S.\ Birman and A.\ Libgober, Contemporary Mathematics {\bf 78},
A.\ M.\ S., Providence, Rhode Island, 1988, pp.\ 705-730.

\end{thebibliography}
}
\end{document}